\documentclass[10pt]{article}

\font\smallit=cmti10
\font\smalltt=cmtt10
\font\smallrm=cmr9

\usepackage{amssymb,latexsym,amsmath,epsfig,amsthm}
\usepackage{amsfonts}
\usepackage[justification=centering]{caption}
\usepackage{float}
\usepackage{afterpage}

\makeatletter

\renewcommand\section{\@startsection {section}{1}{\z@}
{-30pt \@plus -1ex \@minus -.2ex}
{2.3ex \@plus.2ex}
{\normalfont\normalsize\bfseries\boldmath}}

\renewcommand\subsection{\@startsection{subsection}{2}{\z@}
{-3.25ex\@plus -1ex \@minus -.2ex}
{1.5ex \@plus .2ex}
{\normalfont\normalsize\bfseries\boldmath}}

\renewcommand{\@seccntformat}[1]{\csname the#1\endcsname. }

\makeatother

\newtheorem{theorem}{Theorem}

\newtheorem{definition}{Definition}

\allowdisplaybreaks

\begin{document}

\begin{center}
{\Large \bf Consecutive primes}
\vskip 5pt
{\Large \bf which are widely digitally delicate}
\vskip 20pt
{\bf Michael Filaseta}\\
{\smallit Dept.~Mathematics, 
University of South Carolina, 
Columbia, SC 29208, USA}\\
{\tt filaseta@math.sc.edu}\\ 
\vskip 20pt
{\bf Jacob Juillerat}\\
{\smallit Dept.~Mathematics, 
University of South Carolina, 
Columbia, SC 29208, USA}\\
{\tt juillerj@email.sc.edu}\\ 
\end{center}

\vskip 5pt
\centerline{\phantom{\smallit Received: , Revised: , Accepted: , Published: }} 

\vskip 12pt 
\centerline{\textit{Dedicated to the fond memory of Ronald Graham}}

\vskip 15pt

\centerline{\bf Abstract}
\vskip 5pt\noindent
We show that for every positive integer $k$, there exist $k$ consecutive primes 
having the property that if any digit of any one of the primes, including any of the infinitely many
leading zero digits, is changed, then that prime becomes composite.

\pagestyle{myheadings} 
\thispagestyle{empty} 
\baselineskip=12.875pt 
\vskip 30pt

\section{Introduction}

In 1978, M.~S.~Klamkin \cite{klamkin} posed the following problem.
\vskip 5pt
\centerline{\parbox[t]{12cm}{\textit{Does there exist any prime number such that if any digit (in base $10$) is changed to any other digit, the resulting number is always composite?
}}}
\vskip 8pt\noindent
In addition to computations establishing the existence of such a prime, the published solutions in 1979 to this problem included a proof by P.~Erd\H{o}s \cite{Erd79} that there exist infinitely many such primes. 
Borrowing the terminology from J.~Hopper and P.~Pollack  \cite{hopperpollack}, we call such primes \textit{digitally delicate}.  
The first digitally delicate prime is $294001$. 
Thus, $294001$ is a prime and, for every $d \in \{ 0, 1, \ldots, 9 \}$, each of the numbers
\[
d\hspace{.1em}94001, \quad 2d4001, \quad 29d\hspace{.1em}001,  \quad 294d\hspace{.1em}01, \quad 2940d1, \quad 29400d
\]
is either equal to $294001$ or composite.   The proof provided by Erd\H{o}s consisted of creating a partial covering system of the integers (defined in the next section) followed by a sieve argument. 
In 2011, T.~Tao \cite{tao} showed by refining the sieve argument of Erd\H{o}s that a positive proportion (in terms of asymptotic density) of the primes are digitally delicate.  
In 2013, S.~Konyagin \cite{konyagin}
pointed out that a similar approach implies that a positive proportion of composite numbers $n$, coprime to $10$, satisfy the property that if any digit in the base $10$ 
representation of $n$ is changed, then the resulting number remains composite.   
For example, the number $n=212159$ satisfies this property. 
Thus, every number in the set \[\{d12159, 2d2159, 21d159, 212d59, 2121d9, 21215d: d\in\{0,1, 2, \dots,9\}\}\] is composite. 
Later, in 2016, J.~Hopper and P.~Pollack \cite{hopperpollack} resolved a question of Tao's on digitally delicate primes allowing for an arbitrary but fixed number of digit changes to the beginning and end of the prime.
All of these results and their proofs hold for numbers written in an arbitrary base $b$ rather than base $10$, though the proof provided by Erd\H{o}s \cite{Erd79} only addresses the argument in base $10$.

In 2020, the first author and J.~Southwick \cite{filsou} showed that a positive proportion of primes $p$, are \textit{widely digitally delicate}, which they define as having the property that 
if any digit of $p$, \textit{including any one of the infinitely many leading zeros of $p$}, is replaced by any other digit, then the resulting number is composite.  
The proof was specific to base $10$, though they elaborate on other bases for which the analogous argument produces a similar result, including for example base $31$; 
however, it is not even clear whether widely digitally delicate primes exist in every base.
Observe that the first digitally delicate prime, 294001, is not widely digitally delicate since 10294001 is prime. 
It is of some interest to note that even though a positive proportion of the primes are widely digitally delicate, no specific examples of widely digitally delicate primes are known. 
Later in 2020, the authors with J.~Southwick \cite{filjuisou} gave a related argument showing that there are infinitely many (not necessarily a positive proportion) of composite numbers $n$ in base $10$
such that when any digit is inserted in the decimal expansion of $n$, including between two of the infinitely many leading zeros of $n$ and to the right of the units digit of $n$, the number $n$ remains composite
(see also \cite{Fil10}).

In this paper, we show the following.

\begin{theorem}\label{maintheorem}
For every positive integer $k$, there exist $k$ consecutive primes all of which are widely digitally delicate.
\end{theorem}

Let $\mathcal P$ be a set of primes.  It is not difficult to see that 
if $\mathcal P$ has an asymptotic density of $1$ in the set of primes, then there exist $k$ consecutive primes in $\mathcal P$ for each $k \in \mathbb Z^{+}$.
On the other hand, for every $\varepsilon \in (0,1)$, there exists $\mathcal P$ having asymptotic density $1-\varepsilon$ in the set of primes such that 
there do not exist $k$ consecutive primes in $\mathcal P$ for $k$ sufficiently large (more precisely, for $k \ge 1/\varepsilon$). 
Thus, the prior results stated above are not sufficient to establish Theorem~\ref{maintheorem}.   
The main difficulty in using the prior methods to obtain Theorem~\ref{maintheorem} is in the application of sieve techniques in the prior work.  
We want to bypass the use of sieve techniques and instead give complete covering systems to show that there is an arithmetic progression containing infinitely many primes such that every
prime in the arithmetic progression is a widely digitally delicate prime.  This then gives an alternative proof of the result in \cite{filsou}.  After that, the main
driving force behind the proof of Theorem~\ref{maintheorem}, work of D.~Shiu \cite{shiu}, can be applied.  D.~Shiu \cite{shiu} showed that in any arithmetic
progression containing infinitely many primes (that is, $an+b$ with $\gcd(a,b) = 1$ and $a > 0$) there are arbitrarily long sequences of consecutive primes.  
Thus, once we establish through covering systems that such an arithmetic progression exists where every prime in the arithmetic progression is widely
digitally delicate, D.~Shiu's result immediately applies to finish the proof of Theorem~\ref{maintheorem}.   

Our main focus in this paper is on the proof of Theorem~\ref{maintheorem}.
However, in part, this paper is to emphasize that the remarkable work of Shiu \cite{shiu} provides for a nice application to a number of results 
established via covering systems.  One can also take these applications further by looking at the strengthening of Shiu's work by J.~Maynard \cite{maynard}.  
To illustrate the application of Shiu's work in other context, we give some further examples before closing this introduction. 

A Riesel number is a positive odd integer $k$ with the property that $k \cdot 2^{n}-1$ is composite for all positive integers $n$.
A Sierpi\'nski number is a positive odd integer $k$ with the property that $k \cdot 2^{n}+1$ is composite for all nonnegative integers $n$. 
The existence of such $k$ were established in \cite{riesel} and \cite{sierpinski}, respectively, 
though the former is a rather direct consequence of P.~Erd{\H o}s's work in \cite{pe} and the latter is a somewhat less direct application of this same work, 
an observation made by A.~Schinzel (cf.~\cite{FFK}).
A Brier number is a number $k$ which is simultaneously Riesel and Sierpi\'nski, named after Eric Brier who first considered them (cf.~\cite{FFK}).  
The smallest known Brier number, discovered by Christophe Clavier in 2014 (see \cite{sloantwo}) is
\[
3316923598096294713661.
\]
As is common with all these numbers, examples typically come from covering systems giving an arithmetic progression of examples.  
In particular, Clavier established that every number in the arithmetic progression
\[
3770214739596601257962594704110\,n + 
3316923598096294713661, \quad n \in \mathbb Z^{+} \cup \{ 0 \}
\]
is a Brier number.  Since the numbers $3770214739596601257962594704110$ and $3316923598096294713661$ are coprime, Shiu's theorem gives the following.

\begin{theorem}
For every positive integer $k$, there exist $k$ consecutive primes all of which are Brier numbers.
\end{theorem}

Observe that as an immediate consequence the same result holds if Brier numbers are replaced by Riesel or Sierpi\'nski numbers.  

As another less obvious result to apply Shiu's theorem to, we recall a result of R.~Graham \cite{graham} from 1964.  
He showed that there exist relatively prime positive integers $a$ and $b$ such that the recursive Fibonacci-like sequence
\begin{equation}\label{grahamrecursion}
u_{0} = a, \quad u_{1} = b, \quad \text{and} \quad u_{n+1} = u_{n} + u_{n-1} \quad \text{for integers $n \ge 1$},
\end{equation}
consists entirely of composite numbers.  
The known values for admissible $a$ and $b$ have decreased over the years through the work of
others including D.~Knuth \cite{knuth}, J.~W.~Nicol \cite{nicol} and M.~Vsemirnov \cite{vsemirnov}, the latter giving the smallest known such $a$ and $b$
(but notably the same number of digits as the $a$ and $b$ in \cite{nicol}).
The result has also been generalized to other recursions; see
A~Dubickas, A.~Novikas and J.~\v{S}iurys \cite{dns},
D.~Ismailescu, A.~Ko, C.~Lee and J.~Y.~Park \cite{iklp} and
I.~Lunev \cite{lunev}.
As the Graham result concludes with all $u_{i}$ being composite, the initial elements of the sequence, $a$ and $b$, are composite.  
However, there is still a sense in which one can apply Shiu's result.  To be precise, the smallest known example given by Vsemirnov is
done by taking
\[
a = 106276436867 \quad \text{and} \quad b = 35256392432.
\]
With $u_{j}$ defined as above, one can check that each $u_{j}$ is divisible by a prime from the set
\[
\mathcal P = \{ 2, 3, 5, 7, 11, 17, 19, 23, 31, 41, 47, 61, 107, 181, 541, 1103, 2521 \}.
\]
Setting
\[
N = \prod_{p \in \mathcal P} p = 1821895895860356790898731230,
\]
the value of $a$ and $b$ can be replaced by any integers $a$ and $b$ satisfying
\[
a \equiv 106276436867 {\hskip -4pt}\pmod{N} \quad \text{and} \quad
b \equiv 35256392432 {\hskip -4pt}\pmod{N}.
\]
As $\gcd(106276436867,N) = 31$ and $\gcd(35256392432,N) = 2$, these congruences are equivalent to taking
$a = 31 a'$ and $b = 2 b'$ where $a'$ and $b'$ are integers satisfying
\[
a' \equiv 3428272157 {\hskip -4pt}\pmod{58770835350334090028991330}
\]
and
\[
b' \equiv 17628196216 {\hskip -4pt}\pmod{910947947930178395449365615}.
\]
As a direct application of D.~Shiu's result, we have the following.

\begin{theorem}
For every $k \in \mathbb Z^{+}$, there are $k$ consecutive primes $p_{1}, p_{2}, \ldots, p_{k}$
and $k$ consecutive primes $q_{1}, q_{2}, \ldots, q_{k}$ such that for any $i \in \{ 1, 2, \ldots, k \}$, 
the numbers $a = 31 p_{i}$ and $b = 2 q_{i}$ satisfy 
$\gcd(a,b) = 1$ and have the property that 
the $u_{n}$ defined by \eqref{grahamrecursion}
are all composite.
\end{theorem}

\vskip 0pt \noindent
This latter result is not meant to be particularly significant but rather an indication that Shiu's work
does provide information in cases where covering systems are used to form composite numbers.  

Regarding open problems, 
given the recent excellent works surrounding the non-existence of covering systems of particular forms
(cf.~\cite{BBMST, BBMST2, hough, houghnielsen}),
the authors are not convinced that widely digitally delicate primes exist in every base.
Thus, a tantalizing question is whether they exist or whether a positive proportion of the primes
in every base are  widely digitally delicate.  
In the opposite direction, as noted in \cite{filsou}, Carl Pomerance has asked for an unconditional proof that
there exist infinitely many primes which are not digitally delicate or which are not widely digitally delicate.  
For other open problems in this direction, see the end of the introdutcions in \cite{filjuisou} and \cite{filsou}.

\section{The first steps of the argument}

As noted in the introduction, to prove Theorem~\ref{maintheorem}, 
the work of D.~Shiu \cite{shiu} implies that it suffices to obtain an arithmetic progression $An+B$,
with $A$ and $B$ relatively prime positive integers, such that every prime in the arithmetic progression 
is widely digitally delicate.  We will determine such an $A$ and $B$ by finding relatively prime
positive integers $A$ and $B$ satisfying property ($*$) given by
\vskip 5pt
\centerline{($*$){\ }\parbox[t]{11cm}{If $d \in \{ -9, -8, \ldots, -1 \} \cup \{ 1, 2, \ldots, 9 \}$, then each number
in the set 
\[
\mathcal A_{d} = \big\{ An+B + d\cdot 10^{k}: n \in \mathbb Z^{+}, k \in \mathbb Z^{+} \cup \{ 0 \} \big\}
\]
is composite.}}
\vskip 8pt\noindent
As changing a digit of $An+B$, including any one of its infinitely many leading zero digits, corresponds to adding or subtracting
one of the numbers $1, 2, \ldots, 9$ from a digit of $An+B$, we see that relatively prime positive integers $A$ and $B$ 
satisfying property ($*$) also satisfy the property we want, that every prime in $An+B$ is widely digitally delicate.  

To find relatively prime positive integers $A$ and $B$ satisfying property ($*$), we make use of covering systems which we
define as follows.

\begin{definition}
A covering system (or covering) is a finite set of congruences 
\[
x \equiv a_{1} {\hskip -4pt}\pmod{m_{1}}, \quad x \equiv a_{2} {\hskip -4pt}\pmod{m_{2}}, \quad \ldots, \quad x \equiv a_{r} {\hskip -4pt}\pmod{m_{r}},
\]
where $r \in \mathbb Z^{+}$, each $a_{j} \in \mathbb Z$, and each $m_{j} \in \mathbb Z^{+}$, such that every integer
satisfies at least one congruence in the set of congruences.
\end{definition}

\noindent
In other contexts in the literature, further restrictions can be made on the $m_{j}$, so we emphasize here that we
want to allow for $m_{j} = 1$ and for repeated moduli (so that the $m_{j}$ are not necessarily distinct).  There will
be restrictions on the $m_{j}$ that will arise in the covering systems we build due to the approach we are using.  We
will see these as we proceed. 

For each $d \in \{ -9, -8, \ldots, -1 \} \cup \{ 1, 2, \ldots, 9 \}$, we will create a separate covering system to show that the elements
of $\mathcal A_{d}$ in ($*$) are composite.  
Table~\ref{tablenumcong} indicates, for each $d$, the number of different congruences in the covering system corresponding to $d$. 

\begin{table}[!hbt]
\centering
\caption{Number of congruences for each covering}\label{tablenumcong}
\begin{minipage}{3 cm}
\centering
\begin{tabular}{|c|c|}
\hline
$d$ & \# cong. \\ 
\hline \hline
$-9$ & $232$ \\ \hline 
$-8$ & $441$ \\ \hline 
$-7$ & $1$ \\ \hline 
$-6$ & $257$ \\ \hline 
$-5$ & $268$ \\ \hline 
$-4$ & $1$ \\ \hline 
\end{tabular} 
\end{minipage}
\begin{minipage}{3 cm}
\centering
\begin{tabular}{|c|c|}
\hline
$d$ & \# cong. \\ 
\hline \hline
$-3$ & $739$ \\ \hline 
$-2$ & $289$ \\ \hline 
$-1$ & $1$ \\ \hline 
$1$ & $37$ \\ \hline 
$2$ & $1$ \\ \hline 
$3$ & $203$ \\ \hline 
\end{tabular} 
\end{minipage}
\begin{minipage}{3 cm}
\centering
\begin{tabular}{|c|c|}
\hline
$d$ & \# cong. \\ 
\hline \hline
$4$ & $26$ \\ \hline 
$5$ & $1$ \\ \hline 
$6$ & $19$ \\ \hline 
$7$ & $137$ \\ \hline 
$8$ & $1$ \\ \hline 
$9$ & $4$ \\ \hline 
\end{tabular} 
\end{minipage}
\end{table}

The integers we are covering for each $d$ are the exponents $k$ on $10$ in the definition of $\mathcal A_{d}$.
In other words, we will want to view each exponent $k$ as satisfying one of the congruences in our covering system for a given $\mathcal A_{d}$.  
In the end, the values of $A$ and $B$ will be determined by the congruences we choose for the covering systems as well as certain primes that
arise in our method.  

We clarify that the work on digitally delicate primes in prior work mentioned in the introduction used a partial covering of the integers $k$, 
that is a set of congruences where most but not all integers $k$ satisfy at least one of the congruences, 
together with a sieve argument.
The work in \cite{filsou} on widely digitally delicate primes used covering systems for $d \in \{ 1, 2, \ldots, 9 \}$ and the same approach
of partial coverings and sieves for $d \in \{ -9, -8, \ldots, -1 \}$.  
The work in \cite{filjuisou}, like we will use in this paper, made use of covering systems for all $d \in \{ -9, -8, \ldots, -1 \} \cup \{ 1, 2, \ldots, 9 \}$.   
For \cite{filjuisou}, some of the covering systems could be handled rather easily by taking advantage of the fact that we were looking for
composite numbers satisfying a certain property rather than primes.  

Next, we explain more precisely how we create and take advantage of a covering system for a given fixed $d \in \{ -9, -8, \ldots, -1 \} \cup \{ 1, 2, \ldots, 9 \}$.
We begin with a couple illustrative examples.
Table~\ref{tablenumcong} indicates that a number of the $d$ are handled with just one congruence.  
This is accomplished by taking 
\[
A \equiv 0 {\hskip -5pt}\pmod{3} 
\qquad \text{and} \qquad
B \equiv 1 {\hskip -5pt}\pmod{3}.
\]
Observe that each element of $\mathcal A_{d}$ in ($*$) is divisible by $3$ whenever $d \equiv 2 \pmod{3}$.  
Thus, since $A$ and $B$ are positive, as long as we also have $B > 3$, the elements of $\mathcal A_{d}$ for such $d$ are all composite, which is our goal.  
Note the crucial role of the order of $10$ modulo the prime $3$.  The order is $1$, and the covering system for each of these $d$ is simply $k \equiv 0 \pmod{1}$.  
Every integer satisfies this congruence, so it is a covering system. 
The modulus corresponds to the order of $10$ modulo $3$.  Note also that we cannot use the prime $3$ in an analogous way to cover another digit $d$
because the choices for $A$ and $B$, and hence the congruences on $A$ and $B$ above, are to be independent of $d$.  
For example, if $d = 4$, then $An+B + d\cdot 10^{k} \equiv 1 + 4 \equiv 2 \pmod{3}$ and, hence, $An+B + d\cdot 10^{k}$ will not be divisible by $3$.

As a second illustration, we see from Table~\ref{tablenumcong} that we handle the digit $d = 9$ with $4$ congruences.  The congruences for $d = 9$ are
\[
k \equiv 0 {\hskip -5pt}\pmod{2}, \quad \ 
k \equiv 3 {\hskip -5pt}\pmod{4}, \quad \ 
k \equiv 1 {\hskip -5pt}\pmod{8}, \quad \ 
k \equiv 5 {\hskip -5pt}\pmod{8}.
\]
One easily checks that this is a covering system, that is that every integer $k$ satisfies one of these congruences.  To take advantage of this covering system, 
we choose a different prime $p$ for each congruence with $10$ having order modulo $p$ equal to the modulus.  We used
the prime $11$ with $10$ of order $2$,
the prime $101$ with $10$ of order $4$,
the prime $73$ with $10$ of order $8$, and
the prime $137$ with $10$ of order $8$.
We take $A$ divisible by each of these primes.  For ($*$), with $d = 9$, we want $An+B + 9 \cdot 10^{k}$ composite.  
For $k \equiv 0 \pmod{2}$, we accomplish this by taking $B \equiv 2 \pmod{11}$ and $B > 11$ since then $An+B + 9 \cdot 10^{k} \equiv B + 9 \equiv 0 \pmod{11}$.
For $k \equiv 3 \pmod{4}$, we accomplish this by taking $B \equiv 90 \pmod{101}$ and $B > 101$ since then $An+B + 9 \cdot 10^{k} \equiv 90 + 9 \cdot 10^{3} \equiv 9090 \equiv 0 \pmod{101}$. 
Similarly, for $k \equiv 1 \pmod{8}$ and $B \equiv 56 \pmod{73}$, we obtain  $An+B + 9 \cdot 10^{k} \equiv 0 \pmod{73}$; 
and for $k \equiv 5 \pmod{8}$ and $B \equiv 90 \pmod{137}$, we obtain  $An+B + 9 \cdot 10^{k} \equiv 0 \pmod{137}$. 
Thus, taking $B > 137$, we see that ($*$) holds with $d = 9$.  

Of some significance to our explanations later, we note that we could have interchanged the roles of the primes $73$ and $137$ since $10$ has the same order for each of these primes. 
In other words, we could associate $137$ with the congruence $k \equiv 1 \pmod{8}$ above and associate $73$ with the congruence $k \equiv 5 \pmod{8}$.  
Then for $k \equiv 1 \pmod{8}$ and $B \equiv 47 \pmod{137}$, we would have $An+B + 9 \cdot 10^{k} \equiv 0 \pmod{137}$; and 
for $k \equiv 5 \pmod{8}$ and $B \equiv 17 \pmod{73}$, we would have $An+B + 9 \cdot 10^{k} \equiv 0 \pmod{73}$.  
In general, in our construction of widely digitally delicate primes, we want each congruence 
$k \equiv a \pmod{m}$ in a covering system associated with a prime $p$ for which the order of $10$ modulo $p$ is $m$, 
but how we choose the ordering of those primes (which prime goes to which congruence) for a fixed modulus $m$ is irrelevant. 

For each $d \in \{ -9, -8, \ldots, -1 \} \cup \{ 1, 2, \ldots, 9 \}$, 
we determine a covering system of congruences for $k$, 
where each modulus $m$ corresponds to the order of $10$ modulo some prime $p$.  
This imposes a condition on $A$, namely that $A$ is divisible by each of these primes $p$.
Fixing $d$, a congruence from our covering system $k \equiv a \pmod{m}$, and a corresponding prime $p$ with $10$ having order $m$ modulo $p$, 
we determine $B$ such that $A n + B + d\cdot 10^{k} \equiv B + d\cdot 10^{a} \equiv 0 \pmod{p}$.  
Note that the values of $d$, $a$ and $p$ dictate the congruence condition for $B$ modulo $p$.
Each prime $p$ will correspond to a unique congruence condition $B \equiv - d\cdot 10^{a}  \pmod{p}$, 
so the Chinese Remainder Theorem implies the existence of a $B \in \mathbb Z^{+}$ simultaneously satisfying all 
the congruence conditions modulo primes on $B$.
As long as $B$ is large enough, then the condition ($*$) will hold. 

To make sure that there is a prime of the form $An+B$, we will want $\gcd(A,B) = 1$.  
For $k \equiv a \pmod{m}$ and a corresponding prime $p$ as above, 
 we will have $A$ divisible by $p$ and
$B \equiv - d\cdot 10^{a} \pmod{p}$.
Since $d \in \{ -9, -8, \ldots, -1 \} \cup \{ 1, 2, \ldots, 9 \}$, if $p \ge 11$, then we see that $p \nmid B$. 
We will not be using the primes $p \in \{ 2,5 \}$ as $10$ does not have an order modulo these primes. 
We have already seen that we are using the prime $p = 3$ for $d \equiv 2 \pmod{3}$, so this ensures that $3 \nmid B$. 
We will use $p = 7$ for $d \in \{-9, -8, -6, -5, -3, 3, 4 \}$, which then implies $7 \nmid B$.
Therefore, the condition $\gcd(A,B) = 1$ will hold.

Recall that we used the same congruence and corresponding prime in our covering system for each $d \equiv 2 \pmod{3}$.  
There is no obstacle to repeating a congruence for different $d$ if the corresponding prime, having $10$ of order the modulus, is different.  
But in the case of $d \equiv 2 \pmod{3}$, the same prime $3$ was used for different $d$.  To illustrate how we can repeat the use of a prime, 
we return to how we used the prime $p = 11$ above for $d = 9$.  We ended up with $A \equiv 0 \pmod{11}$ and $B \equiv 2 \pmod{11}$.  
In order for us to take advantage of the prime $p = 11$ for $d$, we therefore want 
$A n + B + d\cdot 10^{k} \equiv 2 + d\cdot 10^{k} \equiv 0 \pmod{11}$.  
It is easy to check that this holds for $(d,k) \in \{ (-9,1), (-2,0), (2,1), (9,0) \}$.  
The case $(d,k)  = (9,0)$ is from our example with $d = 9$ above.
The case $(d,k)  = (2,1)$ does not serve a purpose for us as $d = 2$ was covered by our earlier example using the prime $3$ for all $d \equiv 2 \pmod{3}$. 
The cases where $(d,k) \in \{ (-9,1), (-2,0) \}$ are significant, and we make use of congruences modulo $11$ in the covering systems for $d = -9$ and $d = -2$. 
Thus, we are able to repeat the use of some primes for different values of $d$.  
However, this is not the case for most primes we used.   A complete list of the primes which we were able to use for more than one value of $d$ 
is given in Table~\ref{tablerepeatprimes}, together with the list of corresponding $d$'s.  The function $\rho(m,p)$ in this table will be explained in the next section.

\begin{table}[!hbt]
\centering
\caption{Primes used for more than one digit $d$}\label{tablerepeatprimes}
\begin{minipage}{7 cm}
\centering
\begin{tabular}{|c|c|c|}
\hline
prime & $d$'s & $\rho(m,p)$ \\ 
\hline \hline
$3$ & $-7, -4, -1, 2, 5, 8$ & $1$ \\ \hline
$7$ & $-9, -8, -6, -5, -3, 3, 4$ & $1$ \\ \hline
$11$ & $-9, -2, 9$ & $1$ \\ \hline
$13$ & $-9, -3, 3, 4$ & $2$ \\ \hline
$17$ & $-8, -6, -3, -2, 7$ & $1$ \\ \hline
$19$ & $-6, 4$ & $1$ \\ \hline
$23$ & $-9, -8, -6, -3, 3, 7$ & $1$ \\ \hline
$29$ & $-9, -8, -6, 1, 3$ & $1$ \\ \hline
$31$ & $-8, -2, 6$ & $1$ \\ \hline
$37$ & $3, 4$ & $1$ \\ \hline
$43$ & $-8, -3, 1$ & $1$ \\ \hline
$53$ & $-8, -5, 3$ & $1$ \\ \hline
$61$ & $-6, 3, 6$ & $1$ \\ \hline
$67$ & $-9, 7$ & $1$ \\ \hline
$79$ & $-9, -5$ & $2$ \\ \hline
$89$ & $-6, -3, 7$ & $1$ \\ \hline
$103$ & $-9, -8, -3$ & $1$ \\ \hline
\end{tabular} 
\end{minipage}
\begin{minipage}{5.5 cm}
\centering
\vspace{-.43cm}
\begin{tabular}{|c|c|c|}
\hline
prime & $d$'s & $\rho(m,p)$ \\ 
\hline \hline
$199$ & $-6, -3, 7$ & $1$ \\ \hline
$211$ & $-6, 6$ & $1$ \\ \hline
$241$ & $-6, 6$ & $2$ \\ \hline
$331$ & $-8, 7$ & $1$ \\ \hline
$353$ & $-6, 7$ & $1$ \\ \hline
$409$ & $-8, -3$ & $1$ \\ \hline
$449$ & $-9, 7$ & $2$ \\ \hline
$2161$ & $-6, 6$ & $3$ \\ \hline
$3541$ & $-6, 6$ & $1$ \\ \hline
$9091$ & $-6, 6$ & $1$ \\ \hline
$27961$ & $-6, 6$ & $2$ \\ \hline
$1676321$ & $-6, 6$ & $1$ \\ \hline
$3762091$ & $-6, 6$ & $2$ \\ \hline
$4188901$ & $-6, 6$ & $2$ \\ \hline
$39526741$ & $-6, 6$ & $3$ \\ \hline
$5964848081$ & $-6, 6$ & $2$ \\ \hline
\end{tabular} 
\end{minipage}
\end{table}

Recalling that the modulus in a covering system is equal to the order of $10$ modulo a prime $p$, 
the role of primes and the order of $10$ modulo those primes is significant in coming up with covering systems to deduce ($*$). 
A modulus $m$ can be used in a given covering system as many times as there are primes with $10$ of order $m$.  
Thus, for the covering system for $d = 9$, we saw the modulus $8$ being used twice as there are two primes 
with $10$ of order $8$, namely the primes $73$ and $137$.  
One can look at a list of primitive prime factors of $10^{k}-1$ such as in \cite{brill}, but we needed much more
extensive data than what is contained there.  
Our approach uses that the complete list of primes for which $10$ has a given order $m$ is the same as the list
of primes dividing $\Phi_{m}(10)$ and not dividing $m$ where $\Phi_{m}(x)$ is the $m$-th cyclotomic polynomial (cf.~\cite{brill, filjuisou, filsou}).
We used Magma V2.23-1 on a 2017 MacBook Pro to determine different primes dividing $\Phi_{m}(10)$.  
We did not always get a complete factorization but used that if the remaining unfactored part of $\Phi_{m}(10)$ is composite,
relatively prime to the factored part of $\Phi_{m}(10)$ and $m$, and not a prime power, 
then there must be at least two further distinct prime factors of $\Phi_{m}(10)$.
This allowed us then to determine a lower bound on the number of distinct primes of a given order $m$.  
Though we used most of these in our coverings, sometimes we found extra primes that we did not need to use.  

In total, we made use of $673$ different moduli $m$ and $2596$ different primes dividing $\Phi_{m}(10)$ for such $m$. 
Of the  $2596$ different primes, there are $590$ which came from $295$ composite numbers arising from an unfactored part of some $\Phi_{m}(10)$,
and there are $63$ other composite numbers for which only one prime factor of each of the composite numbers was used.
The largest explicit prime (not coming from the $295+63 = 358$ composite numbers) has $1700$ digits, arising from testing what was
initially a large unfactored part of $\Phi_{m}(10)$ for primality and determining it is a prime.
The largest of the $358$ composite numbers has $17234$ digits.  For obvious reasons, we will avoid listing
these primes and composites in this paper, though to help with verification of the results, we are providing the data
from our computations in \cite{Filweb}; more explicit tables can also be found in \cite{juillerat}.

Table~\ref{orderofprimes} in the appendix gives, for each of the $673$ different moduli $m$, 
the detailed information on the number of distinct primes we used with $10$ of order $m$, 
which we denote by $L(m)$.  Thus, $L(m)$ is a lower bound on the total number of distinct primes with $10$ of order $m$.  
Note that $L(m)$ is less than or equal to the number of distinct primes dividing $\Phi_{m}(10)$ but not dividing $m$.

For each $d \in \{ -9, -8, \ldots, -1 \} \cup \{ 1, 2, \ldots, 9 \}$, the goal is to find a covering system so that ($*$) holds.  
We have already given the covering systems we obtained for $d \equiv 2 \pmod{3}$ and for $d = 9$.  
In the next section and the appendix, we elaborate on the covering systems for the remaining $d$.  
We also explain how the reader can verify the data showing these covering systems satisfy the conditions needed for ($*$).

\section{Finishing the argument}

To finish the argument, we need to present a covering system for each value of $d$ in $\{ -9, -8, \ldots, -1 \} \cup \{ 1, 2, \ldots, 9 \}$ 
as described in the previous section.  For the purposes of keeping the presentation of these covering systems manageable, 
for each $m$ listed in Table~\ref{orderofprimes}, 
we take the $L(m)$ primes we found with $10$ of order $m$ and order them in some way.  
Corresponding to the discussion concerning $d = 9$ and the primes $73$ and $137$, 
the particular ordering is not important to us (for example, increasing order would be fine).  Suppose the primes corresponding to $m$ are 
ordered in some way as $p_{1}, p_{2}, \ldots, p_{L(m)}$.   We define $\rho(p_{j},m) = j$.  
Thus, if $p_{j}$ is the $j$-th prime in our ordering of the primes with $10$ of order $m$, we have $\rho(p_{j},m) = j$. 
The particular values we used for $\rho(p_{j},m)$ is not important to the arguments.  So as to make the entries in 
Table~\ref{tablerepeatprimes} correct, the entries for $\rho(p,m)$ indicate the values we used for those primes.  
For example, Table~\ref{orderofprimes} indicates there are $2$ primes of order $6$.  One of them is $7$. 
Table~\ref{tablerepeatprimes} indicates then that $\rho(7,6) = 1$.  Thus, we put $7$ as the first of the $2$ primes
with $10$ of order $6$.  The other prime with $10$ of order $6$ is $13$, and as Table~\ref{tablerepeatprimes}
indicates we set $13$ as the second of the $2$ primes with $10$ of order $6$.

Tables~\ref{covforminus9}-\ref{covfor9} give the covering systems used for each $d \in \{ -9, -8, \ldots, -1 \} \cup \{ 1, 2, \ldots, 9 \}$ with
$d \not\equiv 2 \pmod{3}$.  Rather than indicating the prime, which in some cases has thousands of digits, corresponding to each congruence 
$k \equiv a \pmod{m}$ listed, we simply wrote the value of $\rho(m,p)$.  As $m$ corresponds to the modulus used in the given congruence 
$k \equiv a \pmod{m}$ and the ordering of the primes is not significant to our arguments (any ordering will do), 
this is enough information to confirm the covering arguments.  

That said, the time consuming task of coming up with the $L(m)$ primes to order for each $m$ is nontrivial (at least at this point in time).  
So that this work does not need to be repeated, a complete list of the $L(m)$ primes for each $m$ is given in \cite{Filweb}.  
Further, the tables in the form of lists can be found there as well, with the third column in each case replaced by the prime we used
with $10$ of order the modulus of the congruence in the second column.  
In the way of clarity, recall that the primes were not explicitly computed in the case that the unfactored part of $\Phi_{m}(10)$ was
tested to be composite; instead the composite number is listed in place of both primes in \cite{Filweb}.

For the remainder of this section, we clarify how to verify the information in Tables~\ref{covforminus9}-\ref{covfor9}.  
We address both verification of the covering systems and the information on the primes as listed in \cite{Filweb}.

\subsection{Covering Verification.}

The most direct way to check that a system $\mathcal C$ of congruences
\[
x \equiv a_{1} {\hskip -5pt}\pmod{m_{1}}, \quad x \equiv a_{2} {\hskip -5pt}\pmod{m_{2}}, \quad \ldots, \quad x \equiv a_{s} {\hskip -5pt}\pmod{m_{s}}
\]
is a covering system is to set $\ell = \text{lcm}(m_{1}, m_{2}, \ldots, m_{s})$ and then to check if every integer in the interval $[0,\ell-1]$ satisfies
at least one congruence in $\mathcal C$.  If not, then $\mathcal C$  is not a covering system.  
If on the other hand, every integer in $[0,\ell-1]$ satisfies a congruence in $\mathcal C$, then $\mathcal C$ is a covering system.  To see the latter,
let $n$ be an arbitrary integer, and write $n = \ell q + r$ where $q$ and $r$ are integers with $0 \le r \le \ell-1$.  
Since $r \in [0,\ell-1]$ satisfies some $x \equiv a_{j} \pmod{m_{j}}$ and since $\ell \equiv 0 \pmod{m_{j}}$, we deduce for this same $j$ that
$n = \ell q + r \equiv a_{j} \pmod{m_{j}}$.  

The above is a satisfactory approach if $\ell$ is not too large. 
For the values of $d$ in $ \{ -9, -8, \ldots, -1 \} \cup \{ 1, 2, \ldots, 9 \}$ with $d \not\equiv 2 \pmod{3}$, the least common multiple $\ell$ given by the congruences in 
Tables~\ref{covforminus9}-\ref{covfor9}
are listed in Table~\ref{lcmtable}.
The maximum prime divisor of $\ell$ is also listed in the fourth column of Table~\ref{lcmtable}. 
The value of $\ell$ can exceed $10^{12}$, so we found a more efficient way to test whether one of our systems $\mathcal C$ of congruences,
where $\ell$ is large, is a covering system.

\begin{table}[!hbt]
\centering
\caption{Least common multiple of the moduli for the coverings in each table}\label{lcmtable}
\begin{minipage}{6.6 cm}
\centering
\begin{tabular}{|c|c|c|c|}
\hline
$d$ & Table & $\ell$ & max $p$ \\ 
\hline \hline
$-9$ & $5$ & $14433138720$ & $31$ \\ \hline 
$-8$ & $6$ & $699847948800$ & $17$ \\ \hline 
$-6$ & $7$ & $1045044000$ & $29$ \\ \hline 
$-5$ & $8$ & $56216160$ & $13$ \\ \hline 
$-3$ & $9$ & $1486147703040$ & $19$ \\ \hline 
$-2$ & $10$ & $321253732800$ & $23$ \\ \hline 
 \end{tabular} 
\end{minipage}
\begin{minipage}{5.6 cm}
\centering
\begin{tabular}{|c|c|c|c|}
\hline
$d$ & Table & $\ell$ & max $p$ \\ 
\hline \hline
$1$ & $11$ & $5040$ & $7$ \\ \hline 
$3$ & $12$ & $133333200$ & $37$ \\ \hline 
$4$ & $13$ & $1296$ & $3$ \\ \hline 
$6$ &  $14$ & $360$ & $5$ \\ \hline 
$7$ &$15$ & $18295200$ & $11$ \\ \hline 
$9$ & $16$ & $8$ & $2$ \\ \hline 
\end{tabular} 
\end{minipage}
\end{table}

Suppose $\ell > 10^{6}$ in Table~\ref{lcmtable} and 
the corresponding collection of congruences coming from the table indicated in the second column is $\mathcal C$.  
Let $q$ be the largest prime divisor of $\ell$ as indicated in the fourth column. 
Let $w = 4 \cdot 3 \cdot 5 \cdot q$.  
This choice of $w$ was selected on the basis of some trial and error; 
other choices are certainly reasonable.
We do however want and have that $w$ divides $\ell$.  
Based on the comments above, we would like to know if every integer in the interval $[0,\ell-1]$ satisfies
at least one congruence in $\mathcal C$. 
The basic idea is to take each $u \in [0,w-1]$ and to consider the integers that are congruent to $u$ modulo $w$ in $[0,\ell-1]$.  
One advantage of doing this is that not every congruence in $\mathcal C$ needs to be considered.  
For example, take $d = -3$.  Then Table~\ref{lcmtable} indicates $\ell = 1486147703040$ 
and Table~\ref{tablenumcong} indicates the number of congruences in $\mathcal C$ is $739$.  
From Table~\ref{covforminus3}, the first few of the congruences in $\mathcal C$ are
\[
k \equiv 4 {\hskip -5pt}\pmod{6}, \quad
k \equiv 5 {\hskip -5pt}\pmod{6}, \quad
k \equiv 0 {\hskip -5pt}\pmod{16}, \quad 
k \equiv 11 {\hskip -5pt}\pmod{21}.
\]
Here, $w = 4 \cdot 3 \cdot 5 \cdot 19 = 1140$.  
If we take $u = 0$, then only the third of these congruences can be satisfied by an integer $k$ congruent to $u$ modulo $w$,
as each of the other ones requires $k \not\equiv 0 \pmod{3}$ whereas $k \equiv u \pmod{w}$ requires $k \equiv 0 \pmod{3}$.  
Let $\mathcal C'$ be the congruences in $\mathcal C$ which are consistent with $k \equiv u \pmod{w}$.  
One can determine these congruences by using that there exist integers satisfying both 
$k \equiv a \pmod{m}$ and $k \equiv u \pmod{w}$ if and only if $a \equiv u \pmod{\gcd(m,w)}$. 

Observe that, with $u \in [0,w-1]$ fixed, we would like to know if each integer $v$ of the form
\begin{equation}\label{vequat}
v = w t + u, \quad \text{ with } 0 \le t \le (\ell/w)-1
\end{equation}
satisfies at least one congruence in $\mathcal C'$. 
The main advantage of this approach is that, as we shall now see, not all $\ell/w$ values of $t$ need to be considered.
First, we note that if $\mathcal C'$ is the empty set, then the integers in \eqref{vequat} are not covered and therefore 
$\mathcal C$ is not a covering system.  
Suppose then that $|\mathcal C'| \ge 1$.
Let $\ell'$ denote the least common multiple of the moduli in $\mathcal C'$.  
Let $\delta = \gcd(w,\ell')$. 
We claim that we need only consider $v = w t + u$ where $0 \le t \le (\ell'/\delta)-1$.  
To see this, suppose we know that every $v = w t + u$ with $0 \le t \le (\ell'/\delta)-1$ satisfies one of the congruences in $\mathcal C'$.  
There are integers $q$, $q'$, $r$ and $r'$ satisfying 
$t = \ell' q' + r'$ where $0 \le r' \le \ell'-1$ and $r' = (\ell'/\delta) q + r$, where $0 \le r \le (\ell'/\delta)-1$.  
Then 
\[
v = w t + u = w \ell' q' + w r' + u = w \ell' q' + (w/\delta) \ell' q + w r + u.
\]
The definition of $\delta$ implies that $w/\delta \in \mathbb Z$.  
As each modulus in $\mathcal C'$ divides $\ell'$, we see that $v$ satisfies a congruence in $\mathcal C'$ if and only if
$w r + u$ does.  Here, $w$ and $u$ are fixed and $0 \le r \le (\ell'/\delta)-1$. 
Thus, we see that for each $u \in [0,w-1]$, we can restrict to determining whether $v$ in \eqref{vequat} satisfies a congruence
in $\mathcal C'$ for $0 \le t \le (\ell'/\delta)-1$.  
Returning to the example of $d = -3$, $\ell = 1486147703040$ and $|\mathcal C| = 739$, 
where $w = 1140$ and we considered $u = 0$, one can check that $|\mathcal C'| = 19$, 
$\ell' = 12640320$, $\delta = w$ and $\ell'/\delta = 11088$.  Thus, what started out as ominously checking whether over $10^{12}$ integers
each satisfy at least one of $739$ different congruences is reduced in the case of $u = 0$ to looking at whether $11088$ integers
each satisfy at least one of $19$ different congruences. 
As $u \in [0,w-1]$ varies, the number of computations does as well.  An extreme case for $d = -3$ occurs for $u = 75$, where
we get $\ell'/\delta = 14325696$ and $|\mathcal C'| = 47$.  As $d$ and $u$ vary, though, this computation becomes manageable for 
determining that we have covering systems for each $d$ in $ \{ -9, -8, \ldots, -1 \} \cup \{ 1, 2, \ldots, 9 \}$ with $d \not\equiv 2 \pmod{3}$
and $\ell > 10^{6}$.  On a 2017 MacBook Plus running Maple 2019 with a 2.3 GHz Dual-Core Intel Core i5 processor, the total cpu time for
determining the systems of congruences in Tables~\ref{covforminus9}-\ref{covfor9} are all covering systems took approximately
$2.9$ cpu hours, with almost all of this time spent on the case $d = -3$ which took $2.7$ hours.  The largest value of $\ell'/\delta$
encountered was $\ell'/\delta = 14325696$ which occurred precisely for $d = -3$ and $u \in \{ 75, 303, 531, 759, 987 \}$.

\subsection{Data check.}

The most cumbersome task for us was the determination of the data in Table~\ref{orderofprimes}.
As noted earlier, although the reader can check the data there directly, we have made the list of primes corresponding to each $m$
available through \cite{Filweb}.   With the list of such primes for each $m$, it is still worth indicating how the data can be checked. 
Recall, in particular, the list of primes is not explicit in the case that there was an unfactored part of $\Phi_{m}(10)$.  
In this subsection, we elaborate on what checks should be and were done.
All computations below were done with the MacBook Pro mentioned at the end of the last subsection and using Magma V2.23-1.

For each modulus $m$ used in our constructions (listed in Table~\ref{orderofprimes}), we made a list of primes
$p_{1}, p_{2}, \ldots, p_{s}$, written in increasing order, together with up to two additional primes $q_{1}$ and $q_{2}$, 
included after $p_{s}$ on the list but not written explicitly (as we will discuss).  Each prime came from a factorization or
partial factorization of $\Phi_{m}(10)$.  
The primes $p_{1}, p_{2}, \ldots, p_{s}$ are the distinct primes appearing in the factored part of $\Phi_{m}(10)$, 
and as noted earlier do not include primes dividing $m$.
In some cases, a complete factorization was found for $\Phi_{m}(10)$.  
For such $m$, there are no additional primes $q_{1}$ and $q_{2}$.
If $\Phi_{m}(10)$ had an unfactored part $Q > 1$ (already tested to be composite), then we checked that
$Q$ is relatively prime to $m p_{1} p_{2} \cdots p_{s}$ and that $Q$ is not of the form $N^{k}$ with
$N \in \mathbb Z^{+}$ and $k$ an integer greater than or equal to $2$.  
As this was always the case for the $Q$ tested, we knew each such $Q$ had two distinct prime factors $q_{1}$ and $q_{2}$.
We deduce that there are at least two more primes $q_{j}$, $j \in \{ 1,2 \}$, different from $p_{1}, p_{2}, \ldots, p_{s}$ for which
$10$ has order $m$ modulo $q_{j}$.  As the data only contains the primes used in the covering systems, 
we only included the primes $q_{1}$ and $q_{2}$ that were used.  Thus, despite $Q$ having at least two distinct prime divisors, 
we may have listed anywhere from $0$ to $2$ of them.  The question arises, however, as to how one can list primes that we do not know;
there are primes $q_{1}$ and $q_{2}$ dividing $Q$, but we were unable to (or chose not to) factor $Q$ to determine them explicitly. 
Instead of listing $q_{1}$ and $q_{2}$ then, we opted to list $Q$.  Thus, for each $m$ we associated a list of one of the forms
\[
[p_{1}, p_{2}, \ldots, p_{s}], \quad [p_{1}, p_{2}, \ldots, p_{s},Q], \quad [p_{1}, p_{2}, \ldots, p_{s},Q,Q],
\]
depending on whether $Q$ either did not exist or we used no prime factor of $Q$, we used one prime factor of $Q$, or we used two prime factors of $Q$, respectively.
It is possible that $s=0$; for example, the lists associated with the moduli $2888$ and $2976$ each take the middle form with no $p_{j}$ and one composite number.

For a fixed $m$, given such a list, say from \cite{Filweb}, one merely needs to check:

\begin{itemize}
\setlength\itemsep{-0.25em}
\item
Each element of the list divides $\Phi_{m}(10)$. 
\item
Each element of the list is relatively prime to $m$.
\item
There is at most one composite number, say $Q > 1$, in the list, which may appear at most twice.  The other numbers in the list are distinct primes.
\item
If the composite number $Q$ exists, then $\gcd(Q, p_{1} p_{2} \cdots p_{s}) = 1$. 
\item
If the composite number $Q$ exists twice, then $Q^{1/k} \not\in \mathbb Z^{+}$ for every integer $k \in [2, \log(Q)/\log(2)]$.
\end{itemize}

\noindent
The upper bond in the last item above is simply because $k > \log(Q)/\log(2)$ implies $1 < Q^{1/k} < 2$ and, hence, $Q^{1/k}$ is not an integer.
For each $m$, the value of $L(m)$ in Table~\ref{orderofprimes} is simply the number of elements in the list associated with $m$.

With the data from the tables in the Appendix, also available in \cite{Filweb} with the indicated primes $p_{1}, p_{2}, \ldots, p_{s}, q_{1}, q_{2}$ depending on $m$ as above, 
some further details need to be checked to fully justify the computations.
We verified that whenever $m$ is used as a modulus in a table, it was associated with one of the primes dividing $\Phi_{m}(10)$. 
Furthermore, for any given $d \in \{ -9, -8, \ldots, -1 \} \cup \{ 1, 2, \ldots, 9 \}$, the complete list of primes used as the congruences vary are distinct, 
noting that $q_{1}$ and $q_{2}$, for a given $m$, will be denoted by the same number $Q$ but represent two distinct prime divisors of $Q$.  
As $d$ varies, a given modulus $m$ and a prime $p$ dividing $\Phi_{m}(10)$ can be used more than once as indicated in Table~\ref{tablerepeatprimes}.  
To elaborate, suppose such an $m$ and $p$ is used for each $d \in \mathcal D \subseteq \{ -9, -8, \ldots, -1 \} \cup \{ 1, 2, \ldots, 9 \}$.
For each $d \in \mathcal D$, then, there corresponds a congruence $k \equiv a \pmod{m}$, where $a = a(d)$ will depend on $d$, as well as $m$ and $p$. 
As noted earlier, this is permissible if and only if the values of $d \cdot 10^{a(d)}$ are congruent modulo $p$ for all $d \in \mathcal D$.   
Thus, for each $p$ that occurs in more than one table, as in Table~\ref{tablerepeatprimes}, a check is done to verify the corresponding values of $d \cdot 10^{a(d)}$ 
are congruent modulo $p$.

The verification of the covering systems needed for Theorem~\ref{maintheorem} is complete, and the work of D.~Shiu \cite{shiu} now
implies the theorem.

\vskip 20pt
\centerline{\textbf{\large Appendix}}
\vskip 12pt \noindent
This appendix begins with Table~\ref{orderofprimes}, which gives a lower bound $L(m)$ on 
the number of distinct prime divisors of $\Phi_{m}(10)$ that are relatively prime to $m$.  
The $m$ listed correspond to moduli used in our coverings.  The number $L(m)$ also provides a
lower bound on the number of primes $p$ for which $10$ has order $m$ modulo $p$.

After Table~\ref{orderofprimes}, the remaining Tables~\ref{covforminus9}-\ref{covfor9}
give the congruences $k \equiv a \pmod{m}$ that form the covering systems we obtained for 
$d \in  \{ -9, -8, \ldots, -1 \} \cup \{ 1, 2, \ldots, 9 \}$
with $d \not\equiv 2 \pmod{3}$.  For each congruence, there is an associated prime coming from the
primes listed in Table~\ref{orderofprimes} and that prime is tabulated in the second columns of
Tables~\ref{covforminus9}-\ref{covfor9} (using the notation $\rho(m,p)$ discussed earlier in this paper).

\newpage

\begin{table}[tp]
\footnotesize
\centering
\caption{Number of primes used, $L = L(m)$, with $10$ of order $m$ (Part I)}\label{orderofprimes}
\begin{minipage}{1.9 cm} 
\centering 
\begin{tabular}{|c|c|} 
\hline
$m$ & $L$ \\  \hline
\hline 
$1$ & $1$ \\ \hline 
$2$ & $1$ \\ \hline 
$3$ & $1$ \\ \hline 
$4$ & $1$ \\ \hline 
$5$ & $2$ \\ \hline 
$6$ & $2$ \\ \hline 
$7$ & $2$ \\ \hline 
$8$ & $2$ \\ \hline 
$9$ & $1$ \\ \hline 
$10$ & $1$ \\ \hline 
$11$ & $2$ \\ \hline 
$12$ & $1$ \\ \hline 
$13$ & $3$ \\ \hline 
$14$ & $1$ \\ \hline 
$15$ & $2$ \\ \hline 
$16$ & $2$ \\ \hline 
$17$ & $2$ \\ \hline 
$18$ & $2$ \\ \hline 
$19$ & $1$ \\ \hline 
$20$ & $2$ \\ \hline 
$21$ & $3$ \\ \hline 
$22$ & $3$ \\ \hline 
$23$ & $1$ \\ \hline 
$24$ & $1$ \\ \hline 
$25$ & $3$ \\ \hline 
$26$ & $2$ \\ \hline 
$27$ & $2$ \\ \hline 
$28$ & $3$ \\ \hline 
$29$ & $5$ \\ \hline 
$30$ & $3$ \\ \hline 
$31$ & $3$ \\ \hline 
$32$ & $2$ \\ \hline 
$33$ & $2$ \\ \hline 
$34$ & $3$ \\ \hline 
$35$ & $3$ \\ \hline 
$36$ & $1$ \\ \hline 
$37$ & $3$ \\ \hline 
$38$ & $1$ \\ \hline 
$39$ & $1$ \\ \hline 
$40$ & $2$ \\ \hline 
$42$ & $3$ \\ \hline 
$44$ & $2$ \\ \hline 
$45$ & $2$ \\ \hline 
$46$ & $4$ \\ \hline 
$48$ & $1$ \\ \hline 
$50$ & $3$ \\ \hline 
$51$ & $4$ \\ \hline 
$52$ & $2$ \\ \hline 
$54$ & $2$ \\ \hline 
$55$ & $4$ \\ \hline 
\end{tabular} 
\end{minipage} 
\begin{minipage}{1.9 cm} 
\centering 
\begin{tabular}{|c|c|} 
\hline
$m$ & $L$ \\  \hline
\hline 
$56$ & $2$ \\ \hline 
$57$ & $3$ \\ \hline 
$58$ & $2$ \\ \hline 
$60$ & $3$ \\ \hline 
$62$ & $1$ \\ \hline 
$63$ & $3$ \\ \hline 
$64$ & $4$ \\ \hline 
$65$ & $2$ \\ \hline 
$66$ & $2$ \\ \hline 
$68$ & $3$ \\ \hline 
$69$ & $3$ \\ \hline 
$70$ & $2$ \\ \hline 
$72$ & $3$ \\ \hline 
$74$ & $3$ \\ \hline 
$75$ & $3$ \\ \hline 
$76$ & $2$ \\ \hline 
$77$ & $4$ \\ \hline 
$78$ & $4$ \\ \hline 
$80$ & $2$ \\ \hline 
$81$ & $5$ \\ \hline 
$84$ & $2$ \\ \hline 
$85$ & $3$ \\ \hline 
$87$ & $3$ \\ \hline 
$88$ & $2$ \\ \hline 
$90$ & $2$ \\ \hline 
$91$ & $7$ \\ \hline 
$92$ & $3$ \\ \hline 
$93$ & $1$ \\ \hline 
$95$ & $5$ \\ \hline 
$96$ & $4$ \\ \hline 
$99$ & $4$ \\ \hline 
$100$ & $4$ \\ \hline 
$102$ & $2$ \\ \hline 
$104$ & $2$ \\ \hline 
$105$ & $3$ \\ \hline 
$108$ & $3$ \\ \hline 
$110$ & $4$ \\ \hline 
$111$ & $3$ \\ \hline 
$112$ & $2$ \\ \hline 
$114$ & $2$ \\ \hline 
$115$ & $5$ \\ \hline 
$116$ & $4$ \\ \hline 
$117$ & $4$ \\ \hline 
$119$ & $4$ \\ \hline 
$120$ & $1$ \\ \hline 
$121$ & $4$ \\ \hline 
$124$ & $2$ \\ \hline 
$125$ & $4$ \\ \hline 
$126$ & $2$ \\ \hline 
$130$ & $2$ \\ \hline 
\end{tabular} 
\end{minipage} 
\begin{minipage}{1.9 cm} 
\centering 
\begin{tabular}{|c|c|} 
\hline
$m$ & $L$ \\  \hline
\hline 
$132$ & $3$ \\ \hline 
$133$ & $3$ \\ \hline 
$135$ & $5$ \\ \hline 
$136$ & $2$ \\ \hline 
$138$ & $3$ \\ \hline 
$140$ & $5$ \\ \hline 
$143$ & $4$ \\ \hline 
$144$ & $2$ \\ \hline 
$145$ & $4$ \\ \hline 
$148$ & $6$ \\ \hline 
$152$ & $4$ \\ \hline 
$153$ & $6$ \\ \hline 
$154$ & $6$ \\ \hline 
$155$ & $3$ \\ \hline 
$156$ & $3$ \\ \hline 
$161$ & $5$ \\ \hline 
$162$ & $4$ \\ \hline 
$165$ & $3$ \\ \hline 
$168$ & $3$ \\ \hline 
$169$ & $3$ \\ \hline 
$170$ & $3$ \\ \hline 
$171$ & $3$ \\ \hline 
$174$ & $3$ \\ \hline 
$175$ & $3$ \\ \hline 
$176$ & $2$ \\ \hline 
$180$ & $3$ \\ \hline 
$182$ & $2$ \\ \hline 
$184$ & $2$ \\ \hline 
$185$ & $4$ \\ \hline 
$186$ & $4$ \\ \hline 
$187$ & $3$ \\ \hline 
$190$ & $3$ \\ \hline 
$192$ & $4$ \\ \hline 
$195$ & $3$ \\ \hline 
$198$ & $2$ \\ \hline 
$203$ & $4$ \\ \hline 
$204$ & $6$ \\ \hline 
$207$ & $3$ \\ \hline 
$208$ & $6$ \\ \hline 
$209$ & $5$ \\ \hline 
$210$ & $3$ \\ \hline 
$216$ & $2$ \\ \hline 
$217$ & $3$ \\ \hline 
$220$ & $6$ \\ \hline 
$221$ & $3$ \\ \hline 
$222$ & $4$ \\ \hline 
$228$ & $7$ \\ \hline 
$230$ & $2$ \\ \hline 
$231$ & $6$ \\ \hline 
$232$ & $5$ \\ \hline 
\end{tabular} 
\end{minipage} 
\begin{minipage}{1.9 cm} 
\centering 
\begin{tabular}{|c|c|} 
\hline
$m$ & $L$ \\  \hline
\hline 
$234$ & $3$ \\ \hline 
$238$ & $3$ \\ \hline 
$240$ & $3$ \\ \hline 
$242$ & $5$ \\ \hline 
$247$ & $5$ \\ \hline 
$248$ & $5$ \\ \hline 
$250$ & $2$ \\ \hline 
$252$ & $2$ \\ \hline 
$253$ & $5$ \\ \hline 
$255$ & $4$ \\ \hline 
$259$ & $5$ \\ \hline 
$260$ & $5$ \\ \hline 
$261$ & $8$ \\ \hline 
$264$ & $7$ \\ \hline 
$266$ & $3$ \\ \hline 
$270$ & $4$ \\ \hline 
$272$ & $2$ \\ \hline 
$273$ & $5$ \\ \hline 
$275$ & $5$ \\ \hline 
$276$ & $5$ \\ \hline 
$280$ & $2$ \\ \hline 
$285$ & $3$ \\ \hline 
$286$ & $5$ \\ \hline 
$289$ & $3$ \\ \hline 
$290$ & $6$ \\ \hline 
$296$ & $3$ \\ \hline 
$297$ & $6$ \\ \hline 
$299$ & $3$ \\ \hline 
$304$ & $3$ \\ \hline 
$306$ & $6$ \\ \hline 
$310$ & $4$ \\ \hline 
$312$ & $5$ \\ \hline 
$315$ & $3$ \\ \hline 
$319$ & $4$ \\ \hline 
$322$ & $6$ \\ \hline 
$323$ & $3$ \\ \hline 
$330$ & $3$ \\ \hline 
$333$ & $4$ \\ \hline 
$336$ & $3$ \\ \hline 
$338$ & $4$ \\ \hline 
$340$ & $5$ \\ \hline 
$341$ & $5$ \\ \hline 
$342$ & $5$ \\ \hline 
$345$ & $6$ \\ \hline 
$348$ & $5$ \\ \hline 
$350$ & $2$ \\ \hline 
$351$ & $4$ \\ \hline 
$352$ & $4$ \\ \hline 
$357$ & $6$ \\ \hline 
$360$ & $3$ \\ \hline 
\end{tabular} 
\end{minipage} 
\begin{minipage}{1.9 cm} 
\centering 
\begin{tabular}{|c|c|} 
\hline
$m$ & $L$ \\  \hline
\hline 
$361$ & $7$ \\ \hline 
$363$ & $3$ \\ \hline 
$364$ & $3$ \\ \hline 
$368$ & $2$ \\ \hline 
$370$ & $2$ \\ \hline 
$372$ & $3$ \\ \hline 
$374$ & $5$ \\ \hline 
$377$ & $4$ \\ \hline 
$380$ & $8$ \\ \hline 
$384$ & $5$ \\ \hline 
$390$ & $3$ \\ \hline 
$391$ & $4$ \\ \hline 
$396$ & $3$ \\ \hline 
$399$ & $5$ \\ \hline 
$403$ & $5$ \\ \hline 
$406$ & $5$ \\ \hline 
$407$ & $6$ \\ \hline 
$408$ & $5$ \\ \hline 
$414$ & $5$ \\ \hline 
$416$ & $5$ \\ \hline 
$418$ & $4$ \\ \hline 
$420$ & $4$ \\ \hline 
$425$ & $3$ \\ \hline 
$429$ & $3$ \\ \hline 
$432$ & $2$ \\ \hline 
$434$ & $5$ \\ \hline 
$435$ & $4$ \\ \hline 
$437$ & $5$ \\ \hline 
$440$ & $2$ \\ \hline 
$442$ & $5$ \\ \hline 
$444$ & $4$ \\ \hline 
$455$ & $5$ \\ \hline 
$456$ & $5$ \\ \hline 
$459$ & $5$ \\ \hline 
$460$ & $8$ \\ \hline 
$462$ & $3$ \\ \hline 
$464$ & $7$ \\ \hline 
$465$ & $5$ \\ \hline 
$476$ & $5$ \\ \hline 
$480$ & $3$ \\ \hline 
$483$ & $7$ \\ \hline 
$484$ & $6$ \\ \hline 
$494$ & $4$ \\ \hline 
$495$ & $6$ \\ \hline 
$496$ & $4$ \\ \hline 
$506$ & $3$ \\ \hline 
$507$ & $7$ \\ \hline 
$510$ & $5$ \\ \hline 
$513$ & $2$ \\ \hline 
$518$ & $2$ \\ \hline 
\end{tabular} 
\end{minipage} 
\begin{minipage}{1.9 cm} 
\centering 
\begin{tabular}{|c|c|} 
\hline
$m$ & $L$ \\  \hline
\hline 
$520$ & $2$ \\ \hline 
$522$ & $4$ \\ \hline 
$527$ & $4$ \\ \hline 
$528$ & $3$ \\ \hline 
$532$ & $5$ \\ \hline 
$540$ & $7$ \\ \hline 
$544$ & $3$ \\ \hline 
$546$ & $2$ \\ \hline 
$552$ & $4$ \\ \hline 
$555$ & $6$ \\ \hline 
$561$ & $4$ \\ \hline 
$570$ & $4$ \\ \hline 
$572$ & $5$ \\ \hline 
$575$ & $6$ \\ \hline 
$578$ & $5$ \\ \hline 
$580$ & $4$ \\ \hline 
$589$ & $2$ \\ \hline 
$592$ & $2$ \\ \hline 
$594$ & $5$ \\ \hline 
$595$ & $4$ \\ \hline 
$598$ & $3$ \\ \hline 
$605$ & $2$ \\ \hline 
$608$ & $3$ \\ \hline 
$609$ & $3$ \\ \hline 
$612$ & $6$ \\ \hline 
$620$ & $2$ \\ \hline 
$621$ & $4$ \\ \hline 
$624$ & $2$ \\ \hline 
$627$ & $4$ \\ \hline 
$630$ & $5$ \\ \hline 
$638$ & $5$ \\ \hline 
$644$ & $3$ \\ \hline 
$646$ & $7$ \\ \hline 
$651$ & $3$ \\ \hline 
$660$ & $8$ \\ \hline 
$663$ & $3$ \\ \hline 
$665$ & $3$ \\ \hline 
$666$ & $4$ \\ \hline 
$672$ & $6$ \\ \hline 
$676$ & $3$ \\ \hline 
$680$ & $4$ \\ \hline 
$682$ & $6$ \\ \hline 
$684$ & $4$ \\ \hline 
$690$ & $5$ \\ \hline 
$693$ & $6$ \\ \hline 
$696$ & $3$ \\ \hline 
$702$ & $4$ \\ \hline 
$704$ & $4$ \\ \hline 
$714$ & $2$ \\ \hline 
$715$ & $3$ \\ \hline
\end{tabular} 
\end{minipage}
\end{table}
\addtocounter{table}{-1}

\begin{table}[tp]
\footnotesize
\centering
\caption{Number of primes used, $L = L(m)$, with $10$ of order $m$ (Part II)}
\begin{minipage}{1.9 cm} 
\centering 
\begin{tabular}{|c|c|} 
\hline
$m$ & $L$ \\  \hline
\hline 
$720$ & $2$ \\ \hline 
$722$ & $3$ \\ \hline 
$726$ & $4$ \\ \hline 
$728$ & $2$ \\ \hline 
$740$ & $10$ \\ \hline 
$741$ & $5$ \\ \hline 
$744$ & $7$ \\ \hline 
$748$ & $4$ \\ \hline 
$754$ & $2$ \\ \hline 
$759$ & $3$ \\ \hline 
$760$ & $2$ \\ \hline 
$765$ & $2$ \\ \hline 
$768$ & $6$ \\ \hline 
$777$ & $5$ \\ \hline 
$782$ & $6$ \\ \hline 
$792$ & $2$ \\ \hline 
$798$ & $3$ \\ \hline 
$805$ & $4$ \\ \hline 
$806$ & $3$ \\ \hline 
$812$ & $6$ \\ \hline 
$814$ & $5$ \\ \hline 
$816$ & $4$ \\ \hline 
$828$ & $2$ \\ \hline 
$833$ & $5$ \\ \hline 
$836$ & $3$ \\ \hline 
$840$ & $5$ \\ \hline 
$845$ & $2$ \\ \hline 
$850$ & $1$ \\ \hline 
$858$ & $5$ \\ \hline 
$867$ & $8$ \\ \hline 
$868$ & $6$ \\ \hline 
$870$ & $7$ \\ \hline 
$874$ & $5$ \\ \hline 
$880$ & $5$ \\ \hline 
$884$ & $4$ \\ \hline 
$888$ & $4$ \\ \hline 
$897$ & $3$ \\ \hline 
$910$ & $5$ \\ \hline 
$912$ & $2$ \\ \hline 
$918$ & $4$ \\ \hline 
$920$ & $2$ \\ \hline 
$924$ & $5$ \\ \hline 
$925$ & $9$ \\ \hline 
$928$ & $4$ \\ \hline 
$930$ & $4$ \\ \hline 
$931$ & $7$ \\ \hline 
$935$ & $5$ \\ \hline 
$952$ & $6$ \\ \hline 
$960$ & $2$ \\ \hline 
$966$ & $3$ \\ \hline 
\end{tabular} 
\end{minipage} 
\begin{minipage}{1.9 cm} 
\centering 
\begin{tabular}{|c|c|} 
\hline
$m$ & $L$ \\  \hline
\hline 
$969$ & $4$ \\ \hline 
$988$ & $3$ \\ \hline 
$990$ & $4$ \\ \hline 
$992$ & $4$ \\ \hline 
$1001$ & $4$ \\ \hline 
$1012$ & $3$ \\ \hline 
$1014$ & $4$ \\ \hline 
$1015$ & $7$ \\ \hline 
$1020$ & $7$ \\ \hline 
$1023$ & $4$ \\ \hline 
$1026$ & $4$ \\ \hline 
$1035$ & $3$ \\ \hline 
$1036$ & $2$ \\ \hline 
$1040$ & $4$ \\ \hline 
$1044$ & $6$ \\ \hline 
$1045$ & $3$ \\ \hline 
$1054$ & $2$ \\ \hline 
$1056$ & $3$ \\ \hline 
$1064$ & $3$ \\ \hline 
$1083$ & $7$ \\ \hline 
$1085$ & $3$ \\ \hline 
$1088$ & $6$ \\ \hline 
$1104$ & $7$ \\ \hline 
$1105$ & $4$ \\ \hline 
$1110$ & $2$ \\ \hline 
$1122$ & $2$ \\ \hline 
$1131$ & $4$ \\ \hline 
$1140$ & $11$ \\ \hline 
$1150$ & $2$ \\ \hline 
$1156$ & $5$ \\ \hline 
$1160$ & $2$ \\ \hline 
$1173$ & $4$ \\ \hline 
$1178$ & $5$ \\ \hline 
$1183$ & $4$ \\ \hline 
$1188$ & $2$ \\ \hline 
$1190$ & $6$ \\ \hline 
$1196$ & $7$ \\ \hline 
$1197$ & $3$ \\ \hline 
$1209$ & $3$ \\ \hline 
$1210$ & $3$ \\ \hline 
$1216$ & $9$ \\ \hline 
$1221$ & $6$ \\ \hline 
$1224$ & $5$ \\ \hline 
$1235$ & $5$ \\ \hline 
$1240$ & $6$ \\ \hline 
$1242$ & $4$ \\ \hline 
$1254$ & $5$ \\ \hline 
$1260$ & $7$ \\ \hline 
$1275$ & $4$ \\ \hline 
$1276$ & $3$ \\ \hline 
\end{tabular} 
\end{minipage} 
\begin{minipage}{1.9 cm} 
\centering 
\begin{tabular}{|c|c|} 
\hline
$m$ & $L$ \\  \hline
\hline 
$1288$ & $7$ \\ \hline 
$1292$ & $3$ \\ \hline 
$1302$ & $4$ \\ \hline 
$1311$ & $3$ \\ \hline 
$1320$ & $8$ \\ \hline 
$1326$ & $2$ \\ \hline 
$1330$ & $3$ \\ \hline 
$1332$ & $6$ \\ \hline 
$1352$ & $8$ \\ \hline 
$1368$ & $4$ \\ \hline 
$1380$ & $9$ \\ \hline 
$1386$ & $3$ \\ \hline 
$1392$ & $3$ \\ \hline 
$1395$ & $5$ \\ \hline 
$1428$ & $2$ \\ \hline 
$1430$ & $4$ \\ \hline 
$1440$ & $4$ \\ \hline 
$1444$ & $4$ \\ \hline 
$1445$ & $6$ \\ \hline 
$1452$ & $1$ \\ \hline 
$1456$ & $6$ \\ \hline 
$1480$ & $5$ \\ \hline 
$1482$ & $2$ \\ \hline 
$1488$ & $3$ \\ \hline 
$1496$ & $4$ \\ \hline 
$1508$ & $4$ \\ \hline 
$1521$ & $1$ \\ \hline 
$1530$ & $6$ \\ \hline 
$1547$ & $3$ \\ \hline 
$1554$ & $6$ \\ \hline 
$1564$ & $5$ \\ \hline 
$1581$ & $6$ \\ \hline 
$1584$ & $2$ \\ \hline 
$1596$ & $4$ \\ \hline 
$1610$ & $6$ \\ \hline 
$1612$ & $5$ \\ \hline 
$1615$ & $2$ \\ \hline 
$1632$ & $4$ \\ \hline 
$1665$ & $3$ \\ \hline 
$1666$ & $4$ \\ \hline 
$1672$ & $4$ \\ \hline 
$1680$ & $4$ \\ \hline 
$1690$ & $4$ \\ \hline 
$1700$ & $6$ \\ \hline 
$1716$ & $4$ \\ \hline 
$1734$ & $6$ \\ \hline 
$1736$ & $2$ \\ \hline 
$1740$ & $9$ \\ \hline 
$1748$ & $5$ \\ \hline 
$1760$ & $2$ \\ \hline 
\end{tabular} 
\end{minipage} 
\begin{minipage}{1.9 cm} 
\centering 
\begin{tabular}{|c|c|} 
\hline
$m$ & $L$ \\  \hline
\hline 
$1767$ & $4$ \\ \hline 
$1768$ & $4$ \\ \hline 
$1776$ & $4$ \\ \hline 
$1785$ & $2$ \\ \hline 
$1794$ & $6$ \\ \hline 
$1805$ & $5$ \\ \hline 
$1824$ & $5$ \\ \hline 
$1848$ & $4$ \\ \hline 
$1850$ & $5$ \\ \hline 
$1860$ & $8$ \\ \hline 
$1862$ & $1$ \\ \hline 
$1870$ & $2$ \\ \hline 
$1885$ & $3$ \\ \hline 
$1904$ & $5$ \\ \hline 
$1932$ & $4$ \\ \hline 
$1938$ & $5$ \\ \hline 
$1953$ & $6$ \\ \hline 
$1976$ & $4$ \\ \hline 
$1980$ & $8$ \\ \hline 
$1995$ & $4$ \\ \hline 
$2002$ & $5$ \\ \hline 
$2023$ & $3$ \\ \hline 
$2024$ & $2$ \\ \hline 
$2030$ & $4$ \\ \hline 
$2040$ & $5$ \\ \hline 
$2046$ & $3$ \\ \hline 
$2052$ & $1$ \\ \hline 
$2070$ & $6$ \\ \hline 
$2090$ & $3$ \\ \hline 
$2108$ & $1$ \\ \hline 
$2112$ & $5$ \\ \hline 
$2128$ & $4$ \\ \hline 
$2166$ & $5$ \\ \hline 
$2170$ & $5$ \\ \hline 
$2176$ & $2$ \\ \hline 
$2185$ & $5$ \\ \hline 
$2208$ & $7$ \\ \hline 
$2210$ & $3$ \\ \hline 
$2220$ & $5$ \\ \hline 
$2262$ & $2$ \\ \hline 
$2280$ & $4$ \\ \hline 
$2300$ & $4$ \\ \hline 
$2312$ & $5$ \\ \hline 
$2346$ & $4$ \\ \hline 
$2356$ & $3$ \\ \hline 
$2366$ & $4$ \\ \hline 
$2380$ & $4$ \\ \hline 
$2392$ & $2$ \\ \hline 
$2394$ & $5$ \\ \hline 
$2418$ & $5$ \\ \hline 
\end{tabular} 
\end{minipage} 
\begin{minipage}{1.9 cm} 
\centering 
\begin{tabular}{|c|c|} 
\hline
$m$ & $L$ \\  \hline
\hline 
$2420$ & $6$ \\ \hline 
$2432$ & $3$ \\ \hline 
$2442$ & $3$ \\ \hline 
$2448$ & $2$ \\ \hline 
$2484$ & $2$ \\ \hline 
$2508$ & $2$ \\ \hline 
$2535$ & $4$ \\ \hline 
$2550$ & $2$ \\ \hline 
$2565$ & $4$ \\ \hline 
$2576$ & $3$ \\ \hline 
$2584$ & $8$ \\ \hline 
$2601$ & $4$ \\ \hline 
$2604$ & $5$ \\ \hline 
$2622$ & $6$ \\ \hline 
$2652$ & $4$ \\ \hline 
$2660$ & $2$ \\ \hline 
$2664$ & $4$ \\ \hline 
$2704$ & $4$ \\ \hline 
$2736$ & $3$ \\ \hline 
$2775$ & $7$ \\ \hline 
$2790$ & $1$ \\ \hline 
$2793$ & $4$ \\ \hline 
$2888$ & $1$ \\ \hline 
$2890$ & $3$ \\ \hline 
$2904$ & $2$ \\ \hline 
$2907$ & $3$ \\ \hline 
$2912$ & $3$ \\ \hline 
$2960$ & $5$ \\ \hline 
$2964$ & $3$ \\ \hline 
$2976$ & $1$ \\ \hline 
$3003$ & $3$ \\ \hline 
$3042$ & $4$ \\ \hline 
$3060$ & $7$ \\ \hline 
$3094$ & $4$ \\ \hline 
$3108$ & $1$ \\ \hline 
$3128$ & $3$ \\ \hline 
$3135$ & $9$ \\ \hline 
$3162$ & $2$ \\ \hline 
$3179$ & $5$ \\ \hline 
$3192$ & $3$ \\ \hline 
$3230$ & $3$ \\ \hline 
$3249$ & $3$ \\ \hline 
$3255$ & $4$ \\ \hline 
$3264$ & $4$ \\ \hline 
$3315$ & $5$ \\ \hline 
$3330$ & $6$ \\ \hline 
$3344$ & $4$ \\ \hline 
$3380$ & $4$ \\ \hline 
$3420$ & $5$ \\ \hline 
$3432$ & $2$ \\ \hline 
\end{tabular} 
\end{minipage} 
\begin{minipage}{1.9 cm} 
\centering 
\begin{tabular}{|c|c|} 
\hline
$m$ & $L$ \\  \hline
\hline 
$3468$ & $4$ \\ \hline 
$3480$ & $3$ \\ \hline 
$3496$ & $6$ \\ \hline 
$3534$ & $4$ \\ \hline 
$3549$ & $3$ \\ \hline 
$3570$ & $6$ \\ \hline 
$3627$ & $5$ \\ \hline 
$3648$ & $2$ \\ \hline 
$3696$ & $4$ \\ \hline 
$3700$ & $2$ \\ \hline 
$3720$ & $2$ \\ \hline 
$3724$ & $6$ \\ \hline 
$3740$ & $4$ \\ \hline 
$3770$ & $6$ \\ \hline 
$3808$ & $2$ \\ \hline 
$3876$ & $4$ \\ \hline 
$3960$ & $6$ \\ \hline 
$3990$ & $2$ \\ \hline 
$4004$ & $2$ \\ \hline 
$4046$ & $3$ \\ \hline 
$4060$ & $6$ \\ \hline 
$4080$ & $4$ \\ \hline 
$4104$ & $2$ \\ \hline 
$4180$ & $6$ \\ \hline 
$4224$ & $2$ \\ \hline 
$4256$ & $4$ \\ \hline 
$4332$ & $4$ \\ \hline 
$4352$ & $3$ \\ \hline 
$4356$ & $3$ \\ \hline 
$4370$ & $5$ \\ \hline 
$4416$ & $6$ \\ \hline 
$4420$ & $6$ \\ \hline 
$4440$ & $5$ \\ \hline 
$4560$ & $4$ \\ \hline 
$4641$ & $2$ \\ \hline 
$4692$ & $2$ \\ \hline 
$4752$ & $3$ \\ \hline 
$4788$ & $4$ \\ \hline 
$4836$ & $5$ \\ \hline 
$4864$ & $4$ \\ \hline 
$4896$ & $4$ \\ \hline 
$5005$ & $3$ \\ \hline 
$5016$ & $6$ \\ \hline 
$5070$ & $2$ \\ \hline 
$5130$ & $2$ \\ \hline 
$5168$ & $2$ \\ \hline 
$5202$ & $4$ \\ \hline 
$5320$ & $2$ \\ \hline 
$5328$ & $4$ \\ \hline 
$5472$ & $7$ \\ \hline 
\end{tabular} 
\end{minipage}
\end{table}
\addtocounter{table}{-1}

\begin{table}[!tp]
\footnotesize
\centering
\caption{Number of primes used, $L = L(m)$, with $10$ of order $m$ (Part III)}
\begin{minipage}{1.9 cm} 
\centering 
\begin{tabular}{|c|c|} 
\hline
$m$ & $L$ \\  \hline
\hline 
$5544$ & $3$ \\ \hline 
$5550$ & $6$ \\ \hline 
$5586$ & $4$ \\ \hline 
$5776$ & $2$ \\ \hline 
$5780$ & $10$ \\ \hline 
$5808$ & $4$ \\ \hline 
$5814$ & $3$ \\ \hline 
$5928$ & $2$ \\ \hline 
$5985$ & $3$ \\ \hline 
$6006$ & $6$ \\ \hline 
$6069$ & $2$ \\ \hline 
$6188$ & $3$ \\ \hline 
$6256$ & $3$ \\ \hline 
\end{tabular} 
\end{minipage} 
\begin{minipage}{1.9 cm} 
\centering 
\begin{tabular}{|c|c|} 
\hline
$m$ & $L$ \\  \hline
\hline 
$6270$ & $3$ \\ \hline 
$6324$ & $3$ \\ \hline 
$6358$ & $3$ \\ \hline 
$6384$ & $2$ \\ \hline 
$6460$ & $6$ \\ \hline 
$6498$ & $3$ \\ \hline 
$6510$ & $3$ \\ \hline 
$6528$ & $3$ \\ \hline 
$6630$ & $3$ \\ \hline 
$6660$ & $3$ \\ \hline 
$6840$ & $5$ \\ \hline 
$6936$ & $3$ \\ \hline 
$6960$ & $2$ \\ \hline 
\end{tabular} 
\end{minipage} 
\begin{minipage}{1.9 cm} 
\centering 
\begin{tabular}{|c|c|} 
\hline
$m$ & $L$ \\  \hline
\hline 
$7068$ & $5$ \\ \hline 
$7254$ & $2$ \\ \hline 
$7392$ & $4$ \\ \hline 
$7448$ & $7$ \\ \hline 
$7752$ & $5$ \\ \hline 
$7980$ & $10$ \\ \hline 
$8008$ & $4$ \\ \hline 
$8092$ & $3$ \\ \hline 
$8208$ & $4$ \\ \hline 
$8360$ & $3$ \\ \hline 
$8664$ & $3$ \\ \hline 
$8704$ & $2$ \\ \hline 
$8880$ & $3$ \\ \hline 
\end{tabular} 
\end{minipage} 
\begin{minipage}{1.9 cm} 
\centering 
\begin{tabular}{|c|c|} 
\hline
$m$ & $L$ \\  \hline
\hline 
$9405$ & $6$ \\ \hline 
$9537$ & $4$ \\ \hline 
$9576$ & $5$ \\ \hline 
$9792$ & $4$ \\ \hline 
$10010$ & $4$ \\ \hline 
$10032$ & $4$ \\ \hline 
$10140$ & $10$ \\ \hline 
$10260$ & $5$ \\ \hline 
$10336$ & $4$ \\ \hline 
$10944$ & $4$ \\ \hline 
$11100$ & $4$ \\ \hline 
$11172$ & $4$ \\ \hline 
$12512$ & $2$ \\ \hline 
\end{tabular} 
\end{minipage} 
\begin{minipage}{1.9 cm} 
\centering 
\begin{tabular}{|c|c|} 
\hline
$m$ & $L$ \\  \hline
\hline 
$12540$ & $4$ \\ \hline 
$12716$ & $5$ \\ \hline 
$12996$ & $8$ \\ \hline 
$13056$ & $2$ \\ \hline 
$14508$ & $4$ \\ \hline 
$15015$ & $3$ \\ \hline 
$15960$ & $4$ \\ \hline 
$16184$ & $4$ \\ \hline 
$16720$ & $2$ \\ \hline 
$18810$ & $3$ \\ \hline 
$19074$ & $2$ \\ \hline 
$20064$ & $4$ \\ \hline 
$22344$ & $4$ \\ \hline 
\end{tabular} 
\end{minipage} 
\begin{minipage}{1.9 cm} 
\centering 
\vspace{-1.75cm}
\begin{tabular}{|c|c|} 
\hline
$m$ & $L$ \\  \hline
\hline 
$25080$ & $1$ \\ \hline 
$25432$ & $6$ \\ \hline 
$25992$ & $6$ \\ \hline 
$30030$ & $6$ \\ \hline 
$37620$ & $2$ \\ \hline 
$51984$ & $2$ \\ \hline 
$60060$ & $6$ \\ \hline 
$75240$ & $6$ \\ \hline
\end{tabular} 
\end{minipage}
\end{table}


\begin{table}[h!] 
\footnotesize 
\centering 
\caption{Covering information for $d = -9$ (Part I)}
\label{covforminus9}
\begin{minipage}{4.1 cm} 
\centering 
\begin{tabular}{|c|c|} 
\hline
congruence & $p$ \\ \hline
\hline 
$k \equiv 1 \pmod{2}$ & $1$ \\ \hline 
$k \equiv 0 \pmod{6}$ & $2$ \\ \hline 
$k \equiv 4 \pmod{6}$ & $1$ \\ \hline 
$k \equiv 6 \pmod{32}$ & $2$ \\ \hline 
$k \equiv 26 \pmod{28}$ & $1$ \\ \hline 
$k \equiv 16 \pmod{22}$ & $1$ \\ \hline 
$k \equiv 26 \pmod{33}$ & $1$ \\ \hline 
$k \equiv 9 \pmod{13}$ & $2$ \\ \hline 
$k \equiv 0 \pmod{31}$ & $1$ \\ \hline 
$k \equiv 1 \pmod{31}$ & $2$ \\ \hline 
$k \equiv 2 \pmod{31}$ & $3$ \\ \hline 
$k \equiv 34 \pmod{62}$ & $1$ \\ \hline 
$k \equiv 35 \pmod{93}$ & $1$ \\ \hline 
$k \equiv 98 \pmod{186}$ & $1$ \\ \hline 
$k \equiv 68 \pmod{186}$ & $2$ \\ \hline 
$k \equiv 38 \pmod{186}$ & $3$ \\ \hline 
$k \equiv 8 \pmod{186}$ & $4$ \\ \hline 
$k \equiv 40 \pmod{124}$ & $1$ \\ \hline 
$k \equiv 102 \pmod{124}$ & $2$ \\ \hline 
$k \equiv 320 \pmod{372}$ & $1$ \\ \hline 
$k \equiv 134 \pmod{372}$ & $2$ \\ \hline 
$k \equiv 104 \pmod{372}$ & $3$ \\ \hline 
$k \equiv 42 \pmod{248}$ & $1$ \\ \hline 
$k \equiv 136 \pmod{248}$ & $2$ \\ \hline 
$k \equiv 74 \pmod{248}$ & $3$ \\ \hline 
$k \equiv 12 \pmod{248}$ & $4$ \\ \hline 
$k \equiv 168 \pmod{248}$ & $5$ \\ \hline 
$k \equiv 602 \pmod{744}$ & $1$ \\ \hline 
\end{tabular} 
\end{minipage} 
\begin{minipage}{4.1 cm} 
\centering 
\vspace{0.00 cm} 
\begin{tabular}{|c|c|} 
\hline
congruence & $p$ \\ \hline
\hline 
$k \equiv 44 \pmod{744}$ & $2$ \\ \hline 
$k \equiv 200 \pmod{744}$ & $4$ \\ \hline 
$k \equiv 386 \pmod{744}$ & $3$ \\ \hline 
$k \equiv 572 \pmod{744}$ & $5$ \\ \hline 
$k \equiv 728 \pmod{744}$ & $6$ \\ \hline 
$k \equiv 170 \pmod{744}$ & $7$ \\ \hline 
$k \equiv 414 \pmod{496}$ & $1$ \\ \hline 
$k \equiv 446 \pmod{496}$ & $2$ \\ \hline 
$k \equiv 478 \pmod{496}$ & $3$ \\ \hline 
$k \equiv 14 \pmod{496}$ & $4$ \\ \hline 
$k \equiv 662 \pmod{992}$ & $1$ \\ \hline 
$k \equiv 694 \pmod{992}$ & $2$ \\ \hline 
$k \equiv 726 \pmod{992}$ & $3$ \\ \hline 
$k \equiv 758 \pmod{992}$ & $4$ \\ \hline 
$k \equiv 356 \pmod{1488}$ & $1$ \\ \hline 
$k \equiv 1100 \pmod{1488}$ & $2$ \\ \hline 
$k \equiv 542 \pmod{1488}$ & $3$ \\ \hline 
$k \equiv 2774 \pmod{2976}$ & $1$ \\ \hline 
$k \equiv 140 \pmod{155}$ & $1$ \\ \hline 
$k \equiv 16 \pmod{155}$ & $2$ \\ \hline 
$k \equiv 47 \pmod{155}$ & $3$ \\ \hline 
$k \equiv 78 \pmod{310}$ & $1$ \\ \hline 
$k \equiv 264 \pmod{310}$ & $2$ \\ \hline 
$k \equiv 110 \pmod{310}$ & $3$ \\ \hline 
$k \equiv 296 \pmod{310}$ & $4$ \\ \hline 
$k \equiv 17 \pmod{465}$ & $1$ \\ \hline 
$k \equiv 203 \pmod{465}$ & $2$ \\ \hline 
$k \equiv 389 \pmod{465}$ & $3$ \\ \hline 
\end{tabular} 
\end{minipage} 
\begin{minipage}{4.1 cm} 
\centering 
\vspace{0.00 cm} 
\begin{tabular}{|c|c|} 
\hline
congruence & $p$ \\ \hline
\hline 
$k \equiv 80 \pmod{465}$ & $4$ \\ \hline 
$k \equiv 266 \pmod{465}$ & $5$ \\ \hline 
$k \equiv 452 \pmod{930}$ & $1$ \\ \hline 
$k \equiv 638 \pmod{930}$ & $2$ \\ \hline 
$k \equiv 824 \pmod{930}$ & $3$ \\ \hline 
$k \equiv 50 \pmod{930}$ & $4$ \\ \hline 
$k \equiv 236 \pmod{620}$ & $1$ \\ \hline 
$k \equiv 546 \pmod{620}$ & $2$ \\ \hline 
$k \equiv 112 \pmod{1240}$ & $1$ \\ \hline 
$k \equiv 1042 \pmod{1240}$ & $2$ \\ \hline 
$k \equiv 732 \pmod{1240}$ & $4$ \\ \hline 
$k \equiv 422 \pmod{1240}$ & $3$ \\ \hline 
$k \equiv 608 \pmod{1240}$ & $5$ \\ \hline 
$k \equiv 1228 \pmod{1240}$ & $6$ \\ \hline 
$k \equiv 1538 \pmod{1860}$ & $1$ \\ \hline 
$k \equiv 1724 \pmod{1860}$ & $2$ \\ \hline 
$k \equiv 794 \pmod{1860}$ & $4$ \\ \hline 
$k \equiv 20 \pmod{1860}$ & $3$ \\ \hline 
$k \equiv 950 \pmod{1860}$ & $5$ \\ \hline 
$k \equiv 1136 \pmod{1860}$ & $6$ \\ \hline 
$k \equiv 206 \pmod{1860}$ & $7$ \\ \hline 
$k \equiv 392 \pmod{1860}$ & $8$ \\ \hline 
$k \equiv 1322 \pmod{3720}$ & $1$ \\ \hline 
$k \equiv 3182 \pmod{3720}$ & $2$ \\ \hline 
$k \equiv 578 \pmod{1395}$ & $1$ \\ \hline 
$k \equiv 113 \pmod{1395}$ & $2$ \\ \hline 
$k \equiv 1043 \pmod{1395}$ & $3$ \\ \hline 
$k \equiv 299 \pmod{1395}$ & $4$ \\ \hline 
\end{tabular} 
\end{minipage} 
\end{table} 
\addtocounter{table}{-1}

\begin{table}[h!] 
\footnotesize 
\centering 
\caption{Covering information for $d = -9$ (Part II)}
\begin{minipage}{4.1 cm} 
\centering 
\begin{tabular}{|c|c|} 
\hline
congruence & $p$ \\ \hline
\hline 
$k \equiv 1229 \pmod{1395}$ & $5$ \\ \hline 
$k \equiv 764 \pmod{2790}$ & $1$ \\ \hline 
$k \equiv 21 \pmod{217}$ & $1$ \\ \hline 
$k \equiv 176 \pmod{217}$ & $2$ \\ \hline 
$k \equiv 114 \pmod{217}$ & $3$ \\ \hline 
$k \equiv 52 \pmod{434}$ & $1$ \\ \hline 
$k \equiv 424 \pmod{434}$ & $2$ \\
\hline 
\end{tabular} 
\end{minipage} 
\begin{minipage}{4.1 cm} 
\centering 
\vspace{0.00 cm} 
\begin{tabular}{|c|c|} 
\hline
congruence & $p$ \\ \hline
\hline 
$k \equiv 1574 \pmod{1953}$ & $6$ \\ \hline 
$k \equiv 242 \pmod{341}$ & $1$ \\ \hline 
$k \equiv 56 \pmod{341}$ & $2$ \\ \hline 
$k \equiv 211 \pmod{341}$ & $3$ \\ \hline 
$k \equiv 25 \pmod{341}$ & $4$ \\ \hline 
$k \equiv 149 \pmod{341}$ & $5$ \\ \hline 
$k \equiv 304 \pmod{682}$ & $1$ \\
\hline 
\end{tabular} 
\end{minipage} 
\begin{minipage}{4.1 cm} 
\centering 
\vspace{-0.75 cm} 
\begin{tabular}{|c|c|} 
\hline
congruence & $p$ \\ \hline
\hline 
$k \equiv 6290 \pmod{7254}$ & $1$ \\ \hline 
$k \equiv 3314 \pmod{7254}$ & $2$ \\ \hline 
$k \equiv 896 \pmod{14508}$ & $1$ \\ \hline 
$k \equiv 8150 \pmod{14508}$ & $2$ \\ \hline 
$k \equiv 5732 \pmod{14508}$ & $3$ \\ \hline 
$k \equiv 12986 \pmod{14508}$ & $4$ \\ \hline 
$k \equiv 10 \pmod{34}$ & $1$ \\
\hline 
\end{tabular} 
\end{minipage} 
\end{table} 

\begin{table}[h!] 
\footnotesize 
\centering 
\caption{Covering information for $d = -8$ (Part I)}
\label{covforminus8}
\begin{minipage}{4.1 cm} 
\centering 
\begin{tabular}{|c|c|} 
\hline
congruence & $p$ \\ \hline
\hline 
$k \equiv 2 \pmod{13}$ & $1$ \\ \hline 
$k \equiv 3 \pmod{21}$ & $1$ \\ \hline 
$k \equiv 19 \pmod{28}$ & $1$ \\ \hline 
$k \equiv 10 \pmod{22}$ & $1$ \\ \hline 
$k \equiv 9 \pmod{110}$ & $1$ \\ \hline 
$k \equiv 0 \pmod{6}$ & $1$ \\ \hline 
$k \equiv 11 \pmod{15}$ & $1$ \\
\hline 
\end{tabular} 
\end{minipage} 
\begin{minipage}{4.1 cm} 
\centering 
\vspace{0.00 cm} 
\begin{tabular}{|c|c|} 
\hline
congruence & $p$ \\ \hline
\hline 
$k \equiv 1569 \pmod{2176}$ & $1$ \\ \hline 
$k \equiv 481 \pmod{2176}$ & $2$ \\ \hline 
$k \equiv 73 \pmod{1088}$ & $1$ \\ \hline 
$k \equiv 209 \pmod{1088}$ & $2$ \\ \hline 
$k \equiv 345 \pmod{1088}$ & $3$ \\ \hline 
$k \equiv 617 \pmod{1088}$ & $4$ \\ \hline 
$k \equiv 753 \pmod{1088}$ & $5$ \\
\hline 
\end{tabular} 
\end{minipage} 
\begin{minipage}{4.1 cm} 
\centering 
\vspace{0.00 cm} 
\begin{tabular}{|c|c|} 
\hline
congruence & $p$ \\ \hline
\hline 
$k \equiv 1454 \pmod{2040}$ & $5$ \\ \hline 
$k \equiv 111 \pmod{612}$ & $1$ \\ \hline 
$k \equiv 9 \pmod{612}$ & $2$ \\ \hline 
$k \equiv 519 \pmod{612}$ & $3$ \\ \hline 
$k \equiv 417 \pmod{612}$ & $4$ \\ \hline 
$k \equiv 315 \pmod{612}$ & $5$ \\ \hline 
$k \equiv 213 \pmod{612}$ & $6$ \\
\hline 
\end{tabular} 
\end{minipage} 
\end{table} 
\addtocounter{table}{-1}

\begin{table}[h!] 
\footnotesize 
\centering 
\caption{Covering information for $d = -8$ (Part II)}
\begin{minipage}{4.1 cm} 
\centering 
\begin{tabular}{|c|c|} 
\hline
congruence & $p$ \\ \hline
\hline 
$k \equiv 181 \pmod{765}$ & $1$ \\ \hline 
$k \equiv 436 \pmod{765}$ & $2$ \\ \hline 
$k \equiv 1456 \pmod{1530}$ & $1$ \\ \hline 
$k \equiv 691 \pmod{1530}$ & $2$ \\ \hline 
$k \equiv 1252 \pmod{1530}$ & $4$ \\ \hline 
$k \equiv 487 \pmod{1530}$ & $3$ \\ \hline 
$k \equiv 742 \pmod{1530}$ & $5$ \\
\hline 
\end{tabular} 
\end{minipage} 
\begin{minipage}{4.1 cm} 
\centering 
\vspace{0.00 cm} 
\begin{tabular}{|c|c|} 
\hline
congruence & $p$ \\ \hline
\hline 
$k \equiv 693 \pmod{1904}$ & $3$ \\ \hline 
$k \equiv 455 \pmod{1904}$ & $4$ \\ \hline 
$k \equiv 1407 \pmod{1904}$ & $5$ \\ \hline 
$k \equiv 64 \pmod{357}$ & $1$ \\ \hline 
$k \equiv 302 \pmod{357}$ & $2$ \\ \hline 
$k \equiv 268 \pmod{357}$ & $3$ \\ \hline 
$k \equiv 149 \pmod{357}$ & $4$ \\
\hline 
\end{tabular} 
\end{minipage} 
\begin{minipage}{4.1 cm} 
\centering 
\vspace{0.00 cm} 
\begin{tabular}{|c|c|} 
\hline
congruence & $p$ \\ \hline
\hline 
$k \equiv 235 \pmod{374}$ & $4$ \\ \hline 
$k \equiv 65 \pmod{374}$ & $5$ \\ \hline 
$k \equiv 456 \pmod{748}$ & $1$ \\ \hline 
$k \equiv 269 \pmod{748}$ & $2$ \\ \hline 
$k \equiv 82 \pmod{748}$ & $3$ \\ \hline 
$k \equiv 643 \pmod{748}$ & $4$ \\ \hline 
$k \equiv 592 \pmod{1496}$ & $1$ \\
\hline 
\end{tabular} 
\end{minipage} 
\end{table} 
\addtocounter{table}{-1}

\begin{table}[h!] 
\footnotesize 
\centering 
\caption{Covering information for $d = -8$ (Part III)}
\begin{minipage}{4.1 cm} 
\centering 
\begin{tabular}{|c|c|} 
\hline
congruence & $p$ \\ \hline
\hline 
$k \equiv 542 \pmod{2652}$ & $3$ \\ \hline 
$k \equiv 2531 \pmod{2652}$ & $4$ \\ \hline 
$k \equiv 270 \pmod{1105}$ & $1$ \\ \hline 
$k \equiv 712 \pmod{1105}$ & $2$ \\ \hline 
$k \equiv 933 \pmod{1105}$ & $3$ \\ \hline 
$k \equiv 49 \pmod{1105}$ & $4$ \\ \hline 
$k \equiv 2701 \pmod{3315}$ & $1$ \\
\hline 
\end{tabular} 
\end{minipage} 
\begin{minipage}{4.1 cm} 
\centering 
\vspace{0.00 cm} 
\begin{tabular}{|c|c|} 
\hline
congruence & $p$ \\ \hline
\hline 
$k \equiv 1699 \pmod{2312}$ & $5$ \\ \hline 
$k \equiv 696 \pmod{1156}$ & $1$ \\ \hline 
$k \equiv 985 \pmod{1156}$ & $2$ \\ \hline 
$k \equiv 118 \pmod{1156}$ & $3$ \\ \hline 
$k \equiv 407 \pmod{1156}$ & $4$ \\ \hline 
$k \equiv 1121 \pmod{1156}$ & $5$ \\ \hline 
$k \equiv 424 \pmod{867}$ & $1$ \\
\hline 
\end{tabular} 
\end{minipage} 
\begin{minipage}{4.1 cm} 
\centering 
\vspace{0.00 cm} 
\begin{tabular}{|c|c|} 
\hline
congruence & $p$ \\ \hline
\hline 
$k \equiv 543 \pmod{6936}$ & $1$ \\ \hline 
$k \equiv 5167 \pmod{6936}$ & $2$ \\ \hline 
$k \equiv 2855 \pmod{6936}$ & $3$ \\ \hline 
$k \equiv 832 \pmod{3468}$ & $1$ \\ \hline 
$k \equiv 2566 \pmod{3468}$ & $2$ \\ \hline 
$k \equiv 1988 \pmod{3468}$ & $3$ \\ \hline 
$k \equiv 254 \pmod{3468}$ & $4$ \\
\hline 
\end{tabular} 
\end{minipage} 
\end{table} 

\begin{table}[h!] 
\footnotesize 
\centering 
\caption{Covering information for $d = -6$ (Part I)}
\label{covforminus6}
\begin{minipage}{4.1 cm} 
\centering 
\begin{tabular}{|c|c|} 
\hline
congruence & $p$ \\ \hline
\hline 
$k \equiv 3 \pmod{6}$ & $1$ \\ \hline 
$k \equiv 89 \pmod{90}$ & $2$ \\ \hline 
$k \equiv 6 \pmod{16}$ & $1$ \\ \hline 
$k \equiv 14 \pmod{28}$ & $1$ \\ \hline 
$k \equiv 9 \pmod{32}$ & $1$ \\ \hline 
$k \equiv 17 \pmod{18}$ & $1$ \\ \hline 
$k \equiv 6 \pmod{22}$ & $1$ \\
\hline 
\end{tabular} 
\end{minipage} 
\begin{minipage}{4.1 cm} 
\centering 
\vspace{0.00 cm} 
\begin{tabular}{|c|c|} 
\hline
congruence & $p$ \\ \hline
\hline 
$k \equiv 319 \pmod{840}$ & $5$ \\ \hline 
$k \equiv 439 \pmod{1680}$ & $1$ \\ \hline 
$k \equiv 1279 \pmod{1680}$ & $2$ \\ \hline 
$k \equiv 1399 \pmod{1680}$ & $3$ \\ \hline 
$k \equiv 559 \pmod{1680}$ & $4$ \\ \hline 
$k \equiv 49 \pmod{280}$ & $1$ \\ \hline 
$k \equiv 189 \pmod{280}$ & $2$ \\
\hline 
\end{tabular} 
\end{minipage} 
\begin{minipage}{4.1 cm} 
\centering 
\vspace{0.00 cm} 
\begin{tabular}{|c|c|} 
\hline
congruence & $p$ \\ \hline
\hline 
$k \equiv 330 \pmod{464}$ & $6$ \\ \hline 
$k \equiv 446 \pmod{464}$ & $7$ \\ \hline 
$k \equiv 388 \pmod{928}$ & $1$ \\ \hline 
$k \equiv 620 \pmod{928}$ & $2$ \\ \hline 
$k \equiv 852 \pmod{928}$ & $3$ \\ \hline 
$k \equiv 156 \pmod{928}$ & $4$ \\ \hline 
$k \equiv 12 \pmod{87}$ & $1$ \\
\hline 
\end{tabular} 
\end{minipage} 
\end{table} 
\addtocounter{table}{-1}

\begin{table}[ht] 
\footnotesize 
\centering 
\caption{Covering information for $d = -6$ (Part II)}
\begin{minipage}{4.1 cm} 
\centering 
\begin{tabular}{|c|c|} 
\hline
congruence & $p$ \\ \hline
\hline 
$k \equiv 20 \pmod{290}$ & $3$ \\ \hline 
$k \equiv 194 \pmod{290}$ & $4$ \\ \hline 
$k \equiv 50 \pmod{290}$ & $5$ \\ \hline 
$k \equiv 224 \pmod{290}$ & $6$ \\ \hline 
$k \equiv 278 \pmod{580}$ & $1$ \\ \hline 
$k \equiv 18 \pmod{580}$ & $2$ \\ \hline 
$k \equiv 338 \pmod{580}$ & $3$ \\ \hline 
$k \equiv 78 \pmod{580}$ & $4$ \\ \hline 
$k \equiv 978 \pmod{1160}$ & $1$ \\ \hline 
$k \equiv 398 \pmod{1160}$ & $2$ \\ \hline 
$k \equiv 225 \pmod{435}$ & $1$ \\ \hline 
$k \equiv 370 \pmod{435}$ & $2$ \\ \hline 
$k \equiv 80 \pmod{435}$ & $3$ \\ \hline 
$k \equiv 399 \pmod{435}$ & $4$ \\ \hline 
$k \equiv 544 \pmod{870}$ & $1$ \\ \hline 
$k \equiv 254 \pmod{870}$ & $2$ \\ \hline 
$k \equiv 690 \pmod{870}$ & $3$ \\ \hline 
$k \equiv 400 \pmod{870}$ & $4$ \\ \hline 
$k \equiv 110 \pmod{870}$ & $5$ \\ \hline 
$k \equiv 864 \pmod{870}$ & $6$ \\ \hline 
$k \equiv 574 \pmod{870}$ & $7$ \\ \hline 
$k \equiv 2024 \pmod{3480}$ & $1$ \\ \hline 
$k \equiv 1154 \pmod{3480}$ & $2$ \\ \hline 
$k \equiv 284 \pmod{3480}$ & $3$ \\ \hline 
$k \equiv 6374 \pmod{6960}$ & $1$ \\ \hline 
$k \equiv 2894 \pmod{6960}$ & $2$ \\ \hline 
$k \equiv 138 \pmod{1740}$ & $1$ \\ \hline 
$k \equiv 718 \pmod{1740}$ & $2$ \\ \hline 
$k \equiv 1298 \pmod{1740}$ & $3$ \\ \hline 
$k \equiv 1038 \pmod{1740}$ & $4$ \\ \hline 
$k \equiv 1618 \pmod{1740}$ & $5$ \\ \hline 
$k \equiv 458 \pmod{1740}$ & $6$ \\ \hline 
$k \equiv 198 \pmod{1740}$ & $7$ \\ \hline 
$k \equiv 778 \pmod{1740}$ & $8$ \\ \hline 
$k \equiv 1358 \pmod{1740}$ & $9$ \\ \hline 
$k \equiv 140 \pmod{1015}$ & $1$ \\ \hline 
\end{tabular} 
\end{minipage} 
\begin{minipage}{4.1 cm} 
\centering 
\vspace{0.00 cm} 
\begin{tabular}{|c|c|} 
\hline
congruence & $p$ \\ \hline
\hline 
$k \equiv 575 \pmod{1015}$ & $2$ \\ \hline 
$k \equiv 1010 \pmod{1015}$ & $3$ \\ \hline 
$k \equiv 430 \pmod{1015}$ & $4$ \\ \hline 
$k \equiv 865 \pmod{1015}$ & $5$ \\ \hline 
$k \equiv 285 \pmod{1015}$ & $6$ \\ \hline 
$k \equiv 720 \pmod{1015}$ & $7$ \\ \hline 
$k \equiv 1764 \pmod{2030}$ & $1$ \\ \hline 
$k \equiv 1184 \pmod{2030}$ & $2$ \\ \hline 
$k \equiv 604 \pmod{2030}$ & $3$ \\ \hline 
$k \equiv 24 \pmod{2030}$ & $4$ \\ \hline 
$k \equiv 3504 \pmod{4060}$ & $1$ \\ \hline 
$k \equiv 1474 \pmod{4060}$ & $2$ \\ \hline 
$k \equiv 2924 \pmod{4060}$ & $3$ \\ \hline 
$k \equiv 894 \pmod{4060}$ & $4$ \\ \hline 
$k \equiv 2344 \pmod{4060}$ & $5$ \\ \hline 
$k \equiv 314 \pmod{4060}$ & $6$ \\ \hline 
$k \equiv 141 \pmod{203}$ & $1$ \\ \hline 
$k \equiv 170 \pmod{203}$ & $2$ \\ \hline 
$k \equiv 199 \pmod{203}$ & $3$ \\ \hline 
$k \equiv 25 \pmod{203}$ & $4$ \\ \hline 
$k \equiv 54 \pmod{406}$ & $1$ \\ \hline 
$k \equiv 286 \pmod{406}$ & $2$ \\ \hline 
$k \equiv 316 \pmod{406}$ & $3$ \\ \hline 
$k \equiv 142 \pmod{406}$ & $4$ \\ \hline 
$k \equiv 374 \pmod{406}$ & $5$ \\ \hline 
$k \equiv 112 \pmod{812}$ & $1$ \\ \hline 
$k \equiv 84 \pmod{812}$ & $2$ \\ \hline 
$k \equiv 200 \pmod{812}$ & $3$ \\ \hline 
$k \equiv 606 \pmod{812}$ & $4$ \\ \hline 
$k \equiv 432 \pmod{812}$ & $5$ \\ \hline 
$k \equiv 26 \pmod{812}$ & $6$ \\ \hline 
$k \equiv 258 \pmod{609}$ & $1$ \\ \hline 
$k \equiv 55 \pmod{609}$ & $2$ \\ \hline 
$k \equiv 461 \pmod{609}$ & $3$ \\ \hline 
$k \equiv 143 \pmod{319}$ & $1$ \\ \hline 
$k \equiv 56 \pmod{319}$ & $2$ \\ \hline 
\end{tabular} 
\end{minipage} 
\begin{minipage}{4.1 cm} 
\centering 
\vspace{-0.35 cm} 
\begin{tabular}{|c|c|} 
\hline
congruence & $p$ \\ \hline
\hline 
$k \equiv 288 \pmod{319}$ & $3$ \\ \hline 
$k \equiv 201 \pmod{319}$ & $4$ \\ \hline 
$k \equiv 346 \pmod{638}$ & $1$ \\ \hline 
$k \equiv 578 \pmod{638}$ & $2$ \\ \hline 
$k \equiv 172 \pmod{638}$ & $3$ \\ \hline 
$k \equiv 404 \pmod{638}$ & $4$ \\ \hline 
$k \equiv 636 \pmod{638}$ & $5$ \\ \hline 
$k \equiv 114 \pmod{1276}$ & $1$ \\ \hline 
$k \equiv 868 \pmod{1276}$ & $2$ \\ \hline 
$k \equiv 230 \pmod{1276}$ & $3$ \\ \hline 
$k \equiv 260 \pmod{377}$ & $1$ \\ \hline 
$k \equiv 144 \pmod{377}$ & $2$ \\ \hline 
$k \equiv 28 \pmod{377}$ & $3$ \\ \hline 
$k \equiv 289 \pmod{377}$ & $4$ \\ \hline 
$k \equiv 550 \pmod{754}$ & $1$ \\ \hline 
$k \equiv 434 \pmod{754}$ & $2$ \\ \hline 
$k \equiv 1072 \pmod{1508}$ & $3$ \\ \hline 
$k \equiv 318 \pmod{1508}$ & $1$ \\ \hline 
$k \equiv 956 \pmod{1508}$ & $2$ \\ \hline 
$k \equiv 202 \pmod{1508}$ & $4$ \\ \hline 
$k \equiv 840 \pmod{1131}$ & $1$ \\ \hline 
$k \equiv 463 \pmod{1131}$ & $2$ \\ \hline 
$k \equiv 86 \pmod{1131}$ & $3$ \\ \hline 
$k \equiv 1101 \pmod{1131}$ & $4$ \\ \hline 
$k \equiv 724 \pmod{2262}$ & $1$ \\ \hline 
$k \equiv 1478 \pmod{2262}$ & $2$ \\ \hline 
$k \equiv 985 \pmod{1885}$ & $1$ \\ \hline 
$k \equiv 608 \pmod{1885}$ & $2$ \\ \hline 
$k \equiv 1739 \pmod{1885}$ & $3$ \\ \hline 
$k \equiv 2000 \pmod{3770}$ & $1$ \\ \hline 
$k \equiv 3508 \pmod{3770}$ & $3$ \\ \hline 
$k \equiv 2754 \pmod{3770}$ & $2$ \\ \hline 
$k \equiv 1130 \pmod{3770}$ & $4$ \\ \hline 
$k \equiv 2638 \pmod{3770}$ & $5$ \\ \hline 
$k \equiv 1884 \pmod{3770}$ & $6$ \\ \hline 
\end{tabular} 
\end{minipage} 
\end{table} 

\begin{table}[h!] 
\footnotesize 
\centering 
\caption{Covering information for $d = -5$ (Part I)}
\label{covforminus5}
\begin{minipage}{4.1 cm} 
\centering 
\begin{tabular}{|c|c|} 
\hline
congruence & $p$ \\ \hline
\hline 
$k \equiv 0 \pmod{13}$ & $1$ \\ \hline 
$k \equiv 1 \pmod{13}$ & $2$ \\ \hline 
$k \equiv 2 \pmod{13}$ & $3$ \\ \hline 
$k \equiv 16 \pmod{26}$ & $1$ \\ \hline 
$k \equiv 3 \pmod{26}$ & $2$ \\ \hline 
$k \equiv 4 \pmod{52}$ & $1$ \\ \hline 
$k \equiv 17 \pmod{52}$ & $2$ \\ \hline 
$k \equiv 82 \pmod{208}$ & $1$ \\ \hline 
$k \equiv 134 \pmod{208}$ & $2$ \\ \hline 
\end{tabular} 
\end{minipage} 
\begin{minipage}{4.1 cm} 
\centering 
\vspace{0.00 cm} 
\begin{tabular}{|c|c|} 
\hline
congruence & $p$ \\ \hline
\hline 
$k \equiv 186 \pmod{208}$ & $3$ \\ \hline 
$k \equiv 30 \pmod{208}$ & $4$ \\ \hline 
$k \equiv 147 \pmod{208}$ & $5$ \\ \hline 
$k \equiv 199 \pmod{208}$ & $6$ \\ \hline 
$k \equiv 43 \pmod{416}$ & $1$ \\ \hline 
$k \equiv 251 \pmod{416}$ & $2$ \\ \hline 
$k \equiv 303 \pmod{416}$ & $3$ \\ \hline 
$k \equiv 95 \pmod{416}$ & $4$ \\ \hline 
$k \equiv 18 \pmod{39}$ & $1$ \\ \hline 
\end{tabular} 
\end{minipage} 
\begin{minipage}{4.1 cm} 
\centering 
\vspace{0.00 cm} 
\begin{tabular}{|c|c|} 
\hline
congruence & $p$ \\ \hline
\hline 
$k \equiv 70 \pmod{78}$ & $1$ \\ \hline 
$k \equiv 31 \pmod{78}$ & $2$ \\ \hline 
$k \equiv 44 \pmod{78}$ & $3$ \\ \hline 
$k \equiv 5 \pmod{78}$ & $4$ \\ \hline 
$k \equiv 1 \pmod{6}$ & $1$ \\ \hline 
$k \equiv 45 \pmod{117}$ & $1$ \\ \hline 
$k \equiv 84 \pmod{117}$ & $2$ \\ \hline 
$k \equiv 6 \pmod{117}$ & $3$ \\ \hline 
$k \equiv 58 \pmod{234}$ & $1$ \\ \hline 
\end{tabular} 
\end{minipage} 
\end{table} 
\addtocounter{table}{-1}

\begin{table}[ht] 
\footnotesize 
\centering 
\caption{Covering information for $d = -5$ (Part II)}
\begin{minipage}{4.1 cm} 
\centering 
\begin{tabular}{|c|c|} 
\hline
congruence & $p$ \\ \hline
\hline 
$k \equiv 136 \pmod{234}$ & $2$ \\ \hline 
$k \equiv 214 \pmod{234}$ & $3$ \\ \hline 
$k \equiv 110 \pmod{117}$ & $4$ \\ \hline 
$k \equiv 32 \pmod{351}$ & $1$ \\ \hline 
$k \equiv 149 \pmod{351}$ & $2$ \\ \hline 
$k \equiv 266 \pmod{351}$ & $3$ \\ \hline 
$k \equiv 305 \pmod{351}$ & $4$ \\
\hline 
\end{tabular} 
\end{minipage} 
\begin{minipage}{4.1 cm} 
\centering 
\vspace{0.00 cm} 
\begin{tabular}{|c|c|} 
\hline
congruence & $p$ \\ \hline
\hline 
$k \equiv 87 \pmod{91}$ & $7$ \\ \hline 
$k \equiv 49 \pmod{182}$ & $1$ \\ \hline 
$k \equiv 140 \pmod{182}$ & $2$ \\ \hline 
$k \equiv 36 \pmod{364}$ & $1$ \\ \hline 
$k \equiv 127 \pmod{364}$ & $2$ \\ \hline 
$k \equiv 309 \pmod{364}$ & $3$ \\ \hline 
$k \equiv 218 \pmod{728}$ & $1$ \\
\hline 
\end{tabular} 
\end{minipage} 
\begin{minipage}{4.1 cm} 
\centering 
\vspace{0.00 cm} 
\begin{tabular}{|c|c|} 
\hline
congruence & $p$ \\ \hline
\hline 
$k \equiv 1427 \pmod{6006}$ & $2$ \\ \hline 
$k \equiv 2428 \pmod{6006}$ & $3$ \\ \hline 
$k \equiv 3429 \pmod{6006}$ & $4$ \\ \hline 
$k \equiv 4430 \pmod{6006}$ & $5$ \\ \hline 
$k \equiv 5431 \pmod{6006}$ & $6$ \\ \hline 
$k \equiv 790 \pmod{5005}$ & $1$ \\ \hline 
$k \equiv 1791 \pmod{5005}$ & $2$ \\
\hline 
\end{tabular} 
\end{minipage} 
\end{table} 
\addtocounter{table}{-1}

\begin{table}[ht] 
\footnotesize 
\centering 
\caption{Covering information for $d = -5$ (Part III)}
\begin{minipage}{4.1 cm} 
\centering 
\begin{tabular}{|c|c|} 
\hline
congruence & $p$ \\ \hline
\hline 
$k \equiv 206 \pmod{286}$ & $5$ \\ \hline 
$k \equiv 141 \pmod{858}$ & $1$ \\ \hline 
$k \equiv 284 \pmod{858}$ & $2$ \\ \hline 
$k \equiv 427 \pmod{858}$ & $3$ \\ \hline 
$k \equiv 570 \pmod{858}$ & $4$ \\ \hline 
$k \equiv 856 \pmod{858}$ & $5$ \\ \hline 
$k \equiv 713 \pmod{1716}$ & $4$ \\ \hline 
$k \equiv 1571 \pmod{3432}$ & $1$ \\ \hline 
$k \equiv 3287 \pmod{3432}$ & $2$ \\ \hline 
$k \equiv 76 \pmod{715}$ & $1$ \\ \hline 
$k \equiv 362 \pmod{715}$ & $2$ \\ \hline 
$k \equiv 505 \pmod{715}$ & $3$ \\ \hline 
$k \equiv 219 \pmod{1430}$ & $1$ \\ \hline 
$k \equiv 648 \pmod{1430}$ & $2$ \\ \hline 
$k \equiv 934 \pmod{1430}$ & $3$ \\ \hline 
$k \equiv 1363 \pmod{1430}$ & $4$ \\ \hline 
$k \equiv 12 \pmod{169}$ & $1$ \\ \hline 
$k \equiv 25 \pmod{169}$ & $2$ \\ \hline 
$k \equiv 38 \pmod{169}$ & $3$ \\ \hline 
$k \equiv 51 \pmod{338}$ & $1$ \\ \hline 
$k \equiv 64 \pmod{338}$ & $2$ \\ \hline 
$k \equiv 220 \pmod{338}$ & $3$ \\ \hline 
$k \equiv 233 \pmod{338}$ & $4$ \\ \hline 
$k \equiv 77 \pmod{507}$ & $2$ \\ \hline 
$k \equiv 90 \pmod{507}$ & $1$ \\ \hline 
$k \equiv 246 \pmod{507}$ & $3$ \\ \hline 
$k \equiv 259 \pmod{507}$ & $4$ \\ \hline 
$k \equiv 415 \pmod{507}$ & $5$ \\ \hline 
$k \equiv 428 \pmod{507}$ & $6$ \\ \hline 
$k \equiv 441 \pmod{507}$ & $7$ \\ \hline 
$k \equiv 103 \pmod{1014}$ & $1$ \\ \hline 
\end{tabular} 
\end{minipage} 
\begin{minipage}{4.1 cm} 
\centering 
\vspace{0.00 cm} 
\begin{tabular}{|c|c|} 
\hline
congruence & $p$ \\ \hline
\hline 
$k \equiv 272 \pmod{1014}$ & $2$ \\ \hline 
$k \equiv 610 \pmod{1014}$ & $3$ \\ \hline 
$k \equiv 779 \pmod{1014}$ & $4$ \\ \hline 
$k \equiv 116 \pmod{676}$ & $1$ \\ \hline 
$k \equiv 285 \pmod{676}$ & $2$ \\ \hline 
$k \equiv 454 \pmod{676}$ & $3$ \\ \hline 
$k \equiv 623 \pmod{2704}$ & $1$ \\ \hline 
$k \equiv 1299 \pmod{2704}$ & $2$ \\ \hline 
$k \equiv 1975 \pmod{2704}$ & $3$ \\ \hline 
$k \equiv 2651 \pmod{2704}$ & $4$ \\ \hline 
$k \equiv 129 \pmod{1352}$ & $1$ \\ \hline 
$k \equiv 298 \pmod{1352}$ & $2$ \\ \hline 
$k \equiv 467 \pmod{1352}$ & $3$ \\ \hline 
$k \equiv 636 \pmod{1352}$ & $4$ \\ \hline 
$k \equiv 805 \pmod{1352}$ & $5$ \\ \hline 
$k \equiv 974 \pmod{1352}$ & $6$ \\ \hline 
$k \equiv 1143 \pmod{1352}$ & $7$ \\ \hline 
$k \equiv 1312 \pmod{1352}$ & $8$ \\ \hline 
$k \equiv 311 \pmod{845}$ & $1$ \\ \hline 
$k \equiv 480 \pmod{845}$ & $2$ \\ \hline 
$k \equiv 142 \pmod{1690}$ & $1$ \\ \hline 
$k \equiv 818 \pmod{1690}$ & $2$ \\ \hline 
$k \equiv 987 \pmod{1690}$ & $3$ \\ \hline 
$k \equiv 1663 \pmod{1690}$ & $4$ \\ \hline 
$k \equiv 649 \pmod{3380}$ & $1$ \\ \hline 
$k \equiv 1494 \pmod{3380}$ & $2$ \\ \hline 
$k \equiv 2339 \pmod{3380}$ & $3$ \\ \hline 
$k \equiv 3184 \pmod{3380}$ & $4$ \\ \hline 
$k \equiv 155 \pmod{1183}$ & $1$ \\ \hline 
$k \equiv 324 \pmod{1183}$ & $2$ \\ \hline 
$k \equiv 493 \pmod{1183}$ & $3$ \\ \hline 
\end{tabular} 
\end{minipage} 
\begin{minipage}{4.1 cm} 
\centering 
\vspace{-0.75 cm} 
\begin{tabular}{|c|c|} 
\hline
congruence & $p$ \\ \hline
\hline 
$k \equiv 1169 \pmod{1183}$ & $4$ \\ \hline 
$k \equiv 662 \pmod{2366}$ & $1$ \\ \hline 
$k \equiv 831 \pmod{2366}$ & $2$ \\ \hline 
$k \equiv 1845 \pmod{2366}$ & $3$ \\ \hline 
$k \equiv 2014 \pmod{2366}$ & $4$ \\ \hline 
$k \equiv 1000 \pmod{3549}$ & $1$ \\ \hline 
$k \equiv 2183 \pmod{3549}$ & $2$ \\ \hline 
$k \equiv 3366 \pmod{3549}$ & $3$ \\ \hline 
$k \equiv 168 \pmod{2535}$ & $1$ \\ \hline 
$k \equiv 675 \pmod{2535}$ & $2$ \\ \hline 
$k \equiv 1182 \pmod{2535}$ & $3$ \\ \hline 
$k \equiv 2196 \pmod{2535}$ & $4$ \\ \hline 
$k \equiv 1689 \pmod{5070}$ & $1$ \\ \hline 
$k \equiv 4224 \pmod{5070}$ & $2$ \\ \hline 
$k \equiv 506 \pmod{1521}$ & $1$ \\ \hline 
$k \equiv 1013 \pmod{3042}$ & $1$ \\ \hline 
$k \equiv 1520 \pmod{3042}$ & $2$ \\ \hline 
$k \equiv 2534 \pmod{3042}$ & $3$ \\ \hline 
$k \equiv 3041 \pmod{3042}$ & $4$ \\ \hline 
$k \equiv 844 \pmod{10140}$ & $1$ \\ \hline 
$k \equiv 1858 \pmod{10140}$ & $2$ \\ \hline 
$k \equiv 2872 \pmod{10140}$ & $3$ \\ \hline 
$k \equiv 3886 \pmod{10140}$ & $4$ \\ \hline 
$k \equiv 4900 \pmod{10140}$ & $8$ \\ \hline 
$k \equiv 5914 \pmod{10140}$ & $5$ \\ \hline 
$k \equiv 6928 \pmod{10140}$ & $6$ \\ \hline 
$k \equiv 7942 \pmod{10140}$ & $7$ \\ \hline 
$k \equiv 8956 \pmod{10140}$ & $9$ \\ \hline 
$k \equiv 9970 \pmod{10140}$ & $10$ \\ \hline 
\end{tabular} 
\end{minipage} 
\end{table} 

\begin{table}[h!] 
\footnotesize 
\centering 
\caption{Covering information for $d = -3$ (Part I)}
\label{covforminus3}
\begin{minipage}{4.1 cm} 
\centering 
\begin{tabular}{|c|c|} 
\hline
congruence & $p$ \\ \hline
\hline 
$k \equiv 4 \pmod{6}$ & $2$ \\ \hline 
$k \equiv 5 \pmod{6}$ & $1$ \\ \hline 
$k \equiv 0 \pmod{16}$ & $1$ \\ \hline 
$k \equiv 11 \pmod{21}$ & $1$ \\ \hline 
$k \equiv 14 \pmod{22}$ & $1$ \\ \hline 
$k \equiv 40 \pmod{44}$ & $1$ \\ \hline 
$k \equiv 96 \pmod{99}$ & $1$ \\ \hline 
$k \equiv 8 \pmod{64}$ & $1$ \\ \hline 
$k \equiv 24 \pmod{64}$ & $2$ \\ \hline 
$k \equiv 40 \pmod{64}$ & $3$ \\ \hline 
$k \equiv 56 \pmod{64}$ & $4$ \\ \hline 
$k \equiv 25 \pmod{210}$ & $1$ \\ \hline 
$k \equiv 151 \pmod{210}$ & $2$ \\ \hline 
$k \equiv 67 \pmod{210}$ & $3$ \\ \hline 
\end{tabular} 
\end{minipage} 
\begin{minipage}{4.1 cm} 
\centering 
\vspace{0.00 cm} 
\begin{tabular}{|c|c|} 
\hline
congruence & $p$ \\ \hline
\hline 
$k \equiv 193 \pmod{420}$ & $1$ \\ \hline 
$k \equiv 403 \pmod{420}$ & $2$ \\ \hline 
$k \equiv 109 \pmod{420}$ & $3$ \\ \hline 
$k \equiv 319 \pmod{420}$ & $4$ \\ \hline 
$k \equiv 18 \pmod{336}$ & $1$ \\ \hline 
$k \equiv 228 \pmod{336}$ & $2$ \\ \hline 
$k \equiv 102 \pmod{336}$ & $3$ \\ \hline 
$k \equiv 522 \pmod{672}$ & $1$ \\ \hline 
$k \equiv 186 \pmod{672}$ & $2$ \\ \hline 
$k \equiv 396 \pmod{672}$ & $3$ \\ \hline 
$k \equiv 60 \pmod{672}$ & $4$ \\ \hline 
$k \equiv 270 \pmod{672}$ & $5$ \\ \hline 
$k \equiv 606 \pmod{672}$ & $6$ \\ \hline 
$k \equiv 81 \pmod{126}$ & $1$ \\ \hline 
\end{tabular} 
\end{minipage} 
\begin{minipage}{4.1 cm} 
\centering 
\vspace{0.00 cm} 
\begin{tabular}{|c|c|} 
\hline
congruence & $p$ \\ \hline
\hline 
$k \equiv 39 \pmod{126}$ & $2$ \\ \hline 
$k \equiv 249 \pmod{252}$ & $1$ \\ \hline 
$k \equiv 123 \pmod{252}$ & $2$ \\ \hline 
$k \equiv 113 \pmod{176}$ & $1$ \\ \hline 
$k \equiv 69 \pmod{176}$ & $2$ \\ \hline 
$k \equiv 201 \pmod{352}$ & $1$ \\ \hline 
$k \equiv 333 \pmod{352}$ & $2$ \\ \hline 
$k \equiv 25 \pmod{352}$ & $3$ \\ \hline 
$k \equiv 157 \pmod{352}$ & $4$ \\ \hline 
$k \equiv 3 \pmod{704}$ & $1$ \\ \hline 
$k \equiv 531 \pmod{704}$ & $2$ \\ \hline 
$k \equiv 355 \pmod{704}$ & $3$ \\ \hline 
$k \equiv 179 \pmod{704}$ & $4$ \\ \hline 
$k \equiv 267 \pmod{1056}$ & $1$ \\ \hline 
\end{tabular} 
\end{minipage} 
\end{table} 
\addtocounter{table}{-1}

\begin{table}[ht] 
\footnotesize 
\centering 
\caption{Covering information for $d = -3$ (Part II)}
\begin{minipage}{4.1 cm} 
\centering 
\begin{tabular}{|c|c|} 
\hline
congruence & $p$ \\ \hline
\hline 
$k \equiv 619 \pmod{1056}$ & $2$ \\ \hline 
$k \equiv 971 \pmod{1056}$ & $3$ \\ \hline 
$k \equiv 795 \pmod{2112}$ & $1$ \\ \hline 
$k \equiv 91 \pmod{2112}$ & $2$ \\ \hline 
$k \equiv 1499 \pmod{2112}$ & $3$ \\ \hline 
$k \equiv 1851 \pmod{2112}$ & $4$ \\ \hline 
$k \equiv 1147 \pmod{2112}$ & $5$ \\
\hline 
\end{tabular} 
\end{minipage} 
\begin{minipage}{4.1 cm} 
\centering 
\vspace{0.00 cm} 
\begin{tabular}{|c|c|} 
\hline
congruence & $p$ \\ \hline
\hline 
$k \equiv 4268 \pmod{7392}$ & $2$ \\ \hline 
$k \equiv 6116 \pmod{7392}$ & $3$ \\ \hline 
$k \equiv 572 \pmod{7392}$ & $4$ \\ \hline 
$k \equiv 5060 \pmod{5544}$ & $1$ \\ \hline 
$k \equiv 3212 \pmod{5544}$ & $2$ \\ \hline 
$k \equiv 1364 \pmod{5544}$ & $3$ \\ \hline 
$k \equiv 0 \pmod{19}$ & $1$ \\
\hline 
\end{tabular} 
\end{minipage} 
\begin{minipage}{4.1 cm} 
\centering 
\vspace{0.00 cm} 
\begin{tabular}{|c|c|} 
\hline
congruence & $p$ \\ \hline
\hline 
$k \equiv 24 \pmod{57}$ & $3$ \\ \hline 
$k \equiv 79 \pmod{114}$ & $1$ \\ \hline 
$k \equiv 98 \pmod{114}$ & $2$ \\ \hline 
$k \equiv 61 \pmod{228}$ & $1$ \\ \hline 
$k \equiv 175 \pmod{228}$ & $2$ \\ \hline 
$k \equiv 194 \pmod{228}$ & $3$ \\ \hline 
$k \equiv 157 \pmod{228}$ & $4$ \\
\hline 
\end{tabular} 
\end{minipage} 
\end{table} 
\addtocounter{table}{-1}

\begin{table}[ht] 
\footnotesize 
\centering 
\caption{Covering information for $d = -3$ (Part III)}
\begin{minipage}{4.1 cm} 
\centering 
\begin{tabular}{|c|c|} 
\hline
congruence & $p$ \\ \hline
\hline 
$k \equiv 3883 \pmod{5472}$ & $6$ \\ \hline 
$k \equiv 2515 \pmod{5472}$ & $7$ \\ \hline 
$k \equiv 10039 \pmod{10944}$ & $1$ \\ \hline 
$k \equiv 6619 \pmod{10944}$ & $2$ \\ \hline 
$k \equiv 4567 \pmod{10944}$ & $3$ \\ \hline 
$k \equiv 1147 \pmod{10944}$ & $4$ \\ \hline 
$k \equiv 140 \pmod{570}$ & $1$ \\
\hline 
\end{tabular} 
\end{minipage} 
\begin{minipage}{4.1 cm} 
\centering 
\vspace{0.00 cm} 
\begin{tabular}{|c|c|} 
\hline
congruence & $p$ \\ \hline
\hline 
$k \equiv 10 \pmod{532}$ & $1$ \\ \hline 
$k \equiv 390 \pmod{532}$ & $2$ \\ \hline 
$k \equiv 314 \pmod{532}$ & $3$ \\ \hline 
$k \equiv 49 \pmod{532}$ & $4$ \\ \hline 
$k \equiv 182 \pmod{532}$ & $5$ \\ \hline 
$k \equiv 276 \pmod{1064}$ & $1$ \\ \hline 
$k \equiv 124 \pmod{1064}$ & $2$ \\
\hline 
\end{tabular} 
\end{minipage} 
\begin{minipage}{4.1 cm} 
\centering 
\vspace{0.00 cm} 
\begin{tabular}{|c|c|} 
\hline
congruence & $p$ \\ \hline
\hline 
$k \equiv 885 \pmod{1197}$ & $2$ \\ \hline 
$k \equiv 87 \pmod{1197}$ & $3$ \\ \hline 
$k \equiv 2215 \pmod{2394}$ & $1$ \\ \hline 
$k \equiv 1417 \pmod{2394}$ & $2$ \\ \hline 
$k \equiv 619 \pmod{2394}$ & $3$ \\ \hline 
$k \equiv 1550 \pmod{2394}$ & $4$ \\ \hline 
$k \equiv 752 \pmod{2394}$ & $5$ \\
\hline 
\end{tabular} 
\end{minipage} 
\end{table} 
\addtocounter{table}{-1}

\begin{table}[ht] 
\footnotesize 
\centering 
\caption{Covering information for $d = -3$ (Part IV)}
\begin{minipage}{4.1 cm} 
\centering 
\begin{tabular}{|c|c|} 
\hline
congruence & $p$ \\ \hline
\hline 
$k \equiv 69 \pmod{798}$ & $2$ \\ \hline 
$k \equiv 734 \pmod{798}$ & $3$ \\ \hline 
$k \equiv 1190 \pmod{1596}$ & $1$ \\ \hline 
$k \equiv 50 \pmod{1596}$ & $2$ \\ \hline 
$k \equiv 506 \pmod{1596}$ & $3$ \\ \hline 
$k \equiv 962 \pmod{1596}$ & $4$ \\ \hline 
$k \equiv 1988 \pmod{3192}$ & $1$ \\
\hline 
\end{tabular} 
\end{minipage} 
\begin{minipage}{4.1 cm} 
\centering 
\vspace{0.00 cm} 
\begin{tabular}{|c|c|} 
\hline
congruence & $p$ \\ \hline
\hline 
$k \equiv 3129 \pmod{5130}$ & $2$ \\ \hline 
$k \equiv 474 \pmod{540}$ & $1$ \\ \hline 
$k \equiv 510 \pmod{540}$ & $2$ \\ \hline 
$k \equiv 186 \pmod{540}$ & $3$ \\ \hline 
$k \equiv 402 \pmod{540}$ & $4$ \\ \hline 
$k \equiv 78 \pmod{540}$ & $5$ \\ \hline 
$k \equiv 294 \pmod{540}$ & $6$ \\
\hline 
\end{tabular} 
\end{minipage} 
\begin{minipage}{4.1 cm} 
\centering 
\vspace{0.00 cm} 
\begin{tabular}{|c|c|} 
\hline
congruence & $p$ \\ \hline
\hline 
$k \equiv 204 \pmod{209}$ & $5$ \\ \hline 
$k \equiv 337 \pmod{418}$ & $1$ \\ \hline 
$k \equiv 261 \pmod{418}$ & $2$ \\ \hline 
$k \equiv 185 \pmod{418}$ & $3$ \\ \hline 
$k \equiv 109 \pmod{418}$ & $4$ \\ \hline 
$k \equiv 546 \pmod{836}$ & $1$ \\ \hline 
$k \equiv 470 \pmod{836}$ & $2$ \\
\hline 
\end{tabular} 
\end{minipage} 
\end{table} 
\addtocounter{table}{-1}

\begin{table}[ht] 
\footnotesize 
\centering 
\caption{Covering information for $d = -3$ (Part V)}
\begin{minipage}{4.1 cm} 
\centering 
\begin{tabular}{|c|c|} 
\hline
congruence & $p$ \\ \hline
\hline 
$k \equiv 414 \pmod{4180}$ & $4$ \\ \hline 
$k \equiv 2010 \pmod{4180}$ & $5$ \\ \hline 
$k \equiv 2846 \pmod{4180}$ & $6$ \\ \hline 
$k \equiv 1668 \pmod{8360}$ & $1$ \\ \hline 
$k \equiv 6684 \pmod{8360}$ & $2$ \\ \hline 
$k \equiv 4100 \pmod{8360}$ & $3$ \\ \hline 
$k \equiv 756 \pmod{16720}$ & $1$ \\
\hline 
\end{tabular} 
\end{minipage} 
\begin{minipage}{4.1 cm} 
\centering 
\vspace{0.00 cm} 
\begin{tabular}{|c|c|} 
\hline
congruence & $p$ \\ \hline
\hline 
$k \equiv 8964 \pmod{9405}$ & $3$ \\ \hline 
$k \equiv 2694 \pmod{9405}$ & $5$ \\ \hline 
$k \equiv 5829 \pmod{9405}$ & $6$ \\ \hline 
$k \equiv 901 \pmod{990}$ & $1$ \\ \hline 
$k \equiv 571 \pmod{990}$ & $2$ \\ \hline 
$k \equiv 241 \pmod{990}$ & $3$ \\ \hline 
$k \equiv 307 \pmod{990}$ & $4$ \\
\hline 
\end{tabular} 
\end{minipage} 
\begin{minipage}{4.1 cm} 
\centering 
\vspace{0.00 cm} 
\begin{tabular}{|c|c|} 
\hline
congruence & $p$ \\ \hline
\hline 
$k \equiv 1308 \pmod{1976}$ & $4$ \\ \hline 
$k \equiv 282 \pmod{741}$ & $1$ \\ \hline 
$k \equiv 529 \pmod{741}$ & $2$ \\ \hline 
$k \equiv 35 \pmod{741}$ & $3$ \\ \hline 
$k \equiv 738 \pmod{741}$ & $4$ \\ \hline 
$k \equiv 453 \pmod{741}$ & $5$ \\ \hline 
$k \equiv 985 \pmod{1482}$ & $1$ \\
\hline 
\end{tabular} 
\end{minipage} 
\end{table} 
\addtocounter{table}{-1}

\begin{table}[ht] 
\footnotesize 
\centering 
\caption{Covering information for $d = -3$ (Part VI)}
\begin{minipage}{4.1 cm} 
\centering 
\begin{tabular}{|c|c|} 
\hline
congruence & $p$ \\ \hline
\hline 
$k \equiv 435 \pmod{969}$ & $2$ \\ \hline 
$k \equiv 606 \pmod{969}$ & $3$ \\ \hline 
$k \equiv 948 \pmod{969}$ & $4$ \\ \hline 
$k \equiv 2088 \pmod{2907}$ & $1$ \\ \hline 
$k \equiv 1119 \pmod{2907}$ & $2$ \\ \hline 
$k \equiv 150 \pmod{2907}$ & $3$ \\ \hline 
$k \equiv 3228 \pmod{3876}$ & $1$ \\ \hline 
$k \equiv 321 \pmod{3876}$ & $2$ \\ \hline 
$k \equiv 1290 \pmod{3876}$ & $3$ \\ \hline 
$k \equiv 2259 \pmod{3876}$ & $4$ \\ \hline 
$k \equiv 1879 \pmod{1938}$ & $1$ \\ \hline 
$k \equiv 1081 \pmod{1938}$ & $2$ \\ \hline 
$k \equiv 283 \pmod{1938}$ & $3$ \\ \hline 
$k \equiv 625 \pmod{1938}$ & $4$ \\ \hline 
$k \equiv 1765 \pmod{1938}$ & $5$ \\ \hline 
$k \equiv 4843 \pmod{5814}$ & $1$ \\ \hline 
$k \equiv 967 \pmod{5814}$ & $2$ \\ \hline 
$k \equiv 2905 \pmod{5814}$ & $3$ \\ \hline 
$k \equiv 6572 \pmod{7752}$ & $1$ \\ \hline 
$k \equiv 3836 \pmod{7752}$ & $2$ \\ \hline 
$k \equiv 6116 \pmod{7752}$ & $4$ \\ \hline 
$k \equiv 3380 \pmod{7752}$ & $3$ \\ \hline 
$k \equiv 644 \pmod{7752}$ & $5$ \\ \hline 
$k \equiv 815 \pmod{1615}$ & $1$ \\ \hline 
$k \equiv 1461 \pmod{1615}$ & $2$ \\ \hline 
$k \equiv 492 \pmod{3230}$ & $1$ \\ \hline 
$k \equiv 2107 \pmod{3230}$ & $2$ \\ \hline 
$k \equiv 1138 \pmod{3230}$ & $3$ \\ \hline 
$k \equiv 2753 \pmod{6460}$ & $2$ \\ \hline 
$k \equiv 5983 \pmod{6460}$ & $1$ \\ \hline 
$k \equiv 1784 \pmod{6460}$ & $3$ \\ \hline 
$k \equiv 169 \pmod{6460}$ & $4$ \\ \hline 
$k \equiv 5014 \pmod{6460}$ & $5$ \\ \hline 
\end{tabular}%
 
\end{minipage} 
\end{table} 

\begin{table}[ht] 
\footnotesize 
\centering 
\caption{Covering information for $d = -2$ (Part I)}
\label{covforminus2}
\begin{minipage}{4.1 cm} 
\centering 
%
 
\end{minipage} 
\end{table} 
\addtocounter{table}{-1}

\clearpage
\begin{table}[t!] 
\footnotesize 
\centering 
\caption{Covering information for $d = -2$ (Part III)}
\begin{minipage}{4.1 cm} 
\centering 
%
 
\end{minipage} 
\end{table} 

\begin{table}[ht] 
\footnotesize 
\centering 
\caption{Covering information for $d = 1$ (Part I)}
\label{covfor1}
\begin{minipage}{4.1 cm} 
\centering 
%
 
\end{minipage} 
\end{table} 
\addtocounter{table}{-1}

\clearpage

\begin{table}[t!] 
\footnotesize 
\centering 
\caption{Covering information for $d = 1$ (Part II)}
\begin{minipage}{4.1 cm} 
\centering 
%
 
\end{minipage} 
\end{table} 

\begin{table}[b!] 
\footnotesize 
\centering 
\caption{Covering information for $d = 3$ (Part I)}
\label{covfor3}
\begin{minipage}{4.1 cm} 
\centering 
%
 
\end{minipage} 
\end{table} 
\addtocounter{table}{-1}

\clearpage

\begin{table}[t!] 
\footnotesize 
\centering 
\caption{Covering information for $d = 3$ (Part II)\\[-6pt]}
\begin{minipage}{4.1 cm} 
\centering 
%
 
\end{minipage} 
\end{table} 

\begin{table}[h!] 
\vskip -10pt
\footnotesize 
\centering 
\caption{Covering information for $d = 4$\\[-6pt]}
\label{covfor4}
\begin{minipage}{4.1 cm} 
\centering 
%
 
\end{minipage} 
\end{table} 

\begin{table}[b!] 
\vskip -400pt
\footnotesize 
\centering 
\caption{Covering information for $d = 6$\\[-6pt]}
\label{covfor6}
\begin{minipage}{4.1 cm} 
\centering 
%
 
\end{minipage} 
\end{table} 

\clearpage

\begin{table}[hbt!]
\footnotesize 
\centering 
\caption{Covering information for $d = 7$\\[-6pt]}
\label{covfor7}
\begin{minipage}{4.1 cm} 
\centering 
%
 
\end{minipage} 
\end{table} 

\begin{table}[b!]
\vskip -400pt
\footnotesize 
\centering 
\caption{Covering information for $d = 9$\\[-6pt]}
\label{covfor9}
\begin{minipage}{4.1 cm} 
\centering 
%
 
\end{minipage} 
\end{table}

\end{document}